\documentclass[10pt,reqno]{article}

\usepackage[dvips]{graphics}
\usepackage{epsfig}
\usepackage{amsmath,amsfonts,amsthm,amssymb,amscd}
\renewcommand{\thefootnote}{\fnsymbol{footnote}}
\newcommand{\grad}{\nabla}
\newcommand{\pd}[1]{\frac{\partial}{\partial #1}}

\newcommand{\oa}[1]{\vec{#1}}
\newcommand{\ncr}[2]{{#1 \choose #2}}

\newtheorem{thm}{Theorem}[section]

\newtheorem{lem}[thm]{Lemma}

\newtheorem{rek}[thm]{Remark}

\newcommand\be{\begin{equation}}
\newcommand\ee{\end{equation}}
\newcommand\bea{\begin{eqnarray}}
\newcommand\eea{\end{eqnarray}}
\newcommand\bi{\begin{itemize}}
\newcommand\ei{\end{itemize}}
\newcommand\ben{\begin{enumerate}}
\newcommand\een{\end{enumerate}}
\newcommand\bc{\begin{center}}
\newcommand\ec{\end{center}}
\newcommand\ba{\begin{array}}
\newcommand\ea{\end{array}}

\newcommand{\foh}{\frac{1}{2}}  
\newcommand{\twocase}[5]{#1 \begin{cases} #2 & \text{#3}\\ #4
&\text{#5} \end{cases}   }

\newcommand{\E}{\mathbb{E}}

\newcommand{\kijt}{k_{ijt}}

\newcommand{\gl}{\lambda}
\newcommand{\gb}{\beta}
\newcommand{\gbi}{\beta_i}
\newcommand{\gbip}{\beta_{i,p}}
\newcommand{\gbiP}{\beta_{i,P}}
\newcommand{\gbio}{\beta_{i,1}}
\newcommand{\bip}{\beta_{i,p}}

\newcommand{\gep}{\epsilon}
\newcommand{\ga}{\alpha}

\newcommand{\inti}{\mbox{H}_i}

\begin{document}

\[\]

\[\]

\begin{center}

\Large{\textbf{Closed-Form Bayesian Inferences for the Logit
Model\\ via Polynomial Expansions}} \normalsize

\bigskip

\normalfont Steven J. Miller,\ Eric T. Bradlow,\ Kevin
Dayaratna\footnote{Steven J. Miller (Email: sjmiller@math.brown.edu)
is a Tamarkin Assistant Professor of Mathematics at Brown
University.\ \ Eric T. Bradlow (Email: ebradlow@wharton.upenn.edu)
is the K. P. Chao Professor, Professor of Marketing and Statistics,
and Academic Director of the Wharton Small Business Development
Center at the Wharton School of the University of Pennsylvania.\ \
Kevin Dayaratna (Email: kdayarat@rhsmith.umd.edu) is a doctoral
student in the Marketing Department at the University of Maryland.
The authors would like to thank Peter S. Fader and Dan Stone for
comments on this manuscript.}


December, 2005
\end{center}




\begin{center}
\bf{Closed-Form Bayesian Inferences for the Logit Model via
Polynomial Expansions}
\end{center}

\bigskip

\begin{center}
\textbf{Abstract}\end{center}

\bigskip

\noindent Articles in Marketing and choice literatures have
demonstrated the need for incorporating person-level heterogeneity
into behavioral models (e.g., logit models for multiple binary
outcomes as studied here).  However, the logit likelihood extended
with a population distribution of heterogeneity doesn't yield
closed-form inferences, and therefore numerical integration
techniques are relied upon (e.g., MCMC methods).\\

\noindent We present here an alternative, closed-form Bayesian
inferences for the logit model, which we obtain by approximating
the logit likelihood via a polynomial expansion, and then positing
a distribution of heterogeneity from a flexible family that is now
conjugate and integrable. For problems where the response
coefficients are independent, choosing the Gamma distribution
leads to rapidly convergent closed-form expansions; if there are
correlations among the coefficients one can still obtain rapidly
convergent closed-form expansions by positing a distribution of
heterogeneity from a Multivariate Gamma distribution. The solution
then comes from the moment generating function of the Multivariate
Gamma distribution or in general from the multivariate heterogeneity
distribution assumed.\\

\noindent Closed-form Bayesian inferences, derivatives (useful for
elasticity calculations), population distribution parameter
estimates (useful for summarization) and starting values (useful
for complicated algorithms) are hence directly available. Two
simulation studies demonstrate the efficacy of our approach.\\

\bigskip



\clearpage  

\renewcommand{\thefootnote}{\arabic{footnote}}

\setcounter{footnote}{0}

\begin{center}
\textbf{INTRODUCTION}
\end{center}

Whether it's the 20,000+ hits based on a www.google.com search or
the 1000+ hits on www.jstor.org or the hundreds of published
papers in a variety of disciplines from Marketing to Economics
(Hausman and McFadden 1984) to Statistics (Albert and Chib, 1993)
to Transportation (Bierlaire et. al, 1997), the logit model plays
a very prominent role in many literatures as a basis for
probabilistic inferences for binary outcome data. In part, this is
due to the ubiquitous nature of binary outcome data, whether it is
choices to buy in a given product category or not, the choice to
go to a given medical provider or not, and the like; and, in part,
it may be due to the link between random utility theory and the
logit model in which binary choices following the logit model are
the outcome of a rational economic maximization of latent utility
with extreme value distributed errors (McFadden, 1974).

One of the recent advances regarding this class of models, which
has made its use even more widespread, is its ability to
incorporate heterogeneity into the response coefficients,
reflecting the fact that individuals are likely to vary on the
attribute coefficients that influence their choices (Rossi and
Allenby, 1993; McCulloch and Rossi, 1994).  Whether this
heterogeneity is modelled in an hierarchical Bayesian fashion
allowing for complete variation (Gelfand et. al, 1990), in a
latent-class way allowing for discrete segments (Kamakura and
Russell, 1989), or by using a finite mixture approach (Train and
McFadden, 2000), incorporating person-level heterogeneity is now
the ``expected'' rather than the ``exception''.

Unfortunately, the added flexibility that heterogeneity allows
comes with a price -- numerical computation and complexity.  That
is, once one combines the logit choice kernel, a Bernoulli random
variable with logit link function, with a heterogeneity
distribution, closed-form inference is unavailable due to the
non-conjugacy of the product Bernoulli likelihood and the
heterogeneity distribution (prior). Therefore, numerical methods
such as quadrature, simulated maximum likelihood (Revelt and
Train, 1998), or Markov chain Monte Carlo methods (Gelman et. al.,
1995) are commonly employed to integrate over the heterogeneity
distribution and obtain inferences for the parameters that govern
the heterogeneity distribution (the so-called, population level
parameters). For instance, in the case of a Gaussian heterogeneity
distribution this requires the marginal integration of the product
Bernoulli logit likelihood with the Gaussian distribution, to
obtain means, variances, and possibly covariances of the prior.
While faster computing and specialized software has made this
feasible, this research considers an alternative to these
approaches, a ``closed-form'' solution.

That is, in this research we consider a closed-form solution to
the heterogeneous binary logit choice problem that involves
approximating the product Bernoulli logit likelihood via a
polynomial expansion (to any specified accuracy), and then
specifying a rich and flexible class of heterogeneity
distributions for the response coefficients (slopes). If the
response coefficients within individuals are independent, we model
them as arising from the Gamma distribution (albeit we demonstrate
how are results can be obtained for any multivariate distribution)
or, more generally, a mixture of Gamma distributions (McDonald and
Butler, 1990); if the response coefficients are not independent we
model them as arising from a Multivariate Gamma distribution,
which allows correlations among the coefficients. We then
integrate, now possible in closed-form, the approximated logit
model with respect to these families. Once the model is integrated
with respect to the heterogeneity distribution, we then can
either: (a) maximize the marginal likelihood and obtain Maximum
Marginal Likelihood (MML) estimates of the population parameters
and utilize them for conditional inferences (the empirical Bayes
approach: Morrison and Schmittlein, 1981; Morris, 1983;
Schmittlein, Morrison, and Columbo, 1987) or (b) in the case where
the parameters of the heterogeneity distribution are set
informatively based on prior information, historical data,
subjective beliefs, or the like, fully Bayesian inferences are
obtainable.

In this manner, as in Bradlow, Hardie, and Fader (2002)  for the
negative-Binomial distribution, and in Everson and Bradlow (2002)
for the beta-binomial distribution, one can effectively
incorporate prior information and allow shrinkage that Bayesian
models attend to, but can also obtain closed-form inferences
\emph{without} Monte Carlo simulation efforts or quadrature that
can be sensitive to the starting values and/or contain significant
simulation error. We demonstrate the efficacy of our approach
using two simulated studies, therefore supporting its use as an
alternative method. In addition, we also demonstrate that as a
by-product of the method, closed-form derivatives of the marginal
distribution are obtained which are often of interest in that they
inform how the distribution (possibly in particular the mean and
variance) of population effects would change as a function of a
change in the decision inputs (i.e. covariates). These derivatives
are also commonly (and directly) used in the computation of
probability elasticities, thus providing the opportunity for
optimization decisions.

The remainder of this paper is laid out as follows. In Section
\ref{problem} we systematically lay out the problem formulation by
deriving the likelihood, which provides the basis of our
polynomial expansion (an application of a geometric series
expansion), as well as discuss the types of data for which our
method is applicable. In particular, the results presented here
(albeit they are generalizable) are most applicable (as we discuss
in Section \ref{problem}) to product categories for which the
binary response rate is either rare (e.g. durable goods purchases
(Bayus, 1992) and mail catalog responses (Anderson and Simester
(2001)), or those for which the frequency of purchase is high
(e.g. orange juice).

Section \ref{secindhouseexpand} is an in-depth analysis of the
case when the response coefficients are drawn from independent
Gamma distributions. Sections
\ref{sec:gammadistranditsgeneralizations} and
\ref{subsec:seriesexppygamma} contain our key integration results
demonstrating the conjugacy of the approximation to the binary
data likelihood and the Gamma family of distributions (Theorem
\ref{thm:univargamma}). Details of the integration lemma and plots
of the robustness of the Gamma family and its generalizations that
we consider are in Appendix \ref{ggamma}. We discuss computational
issues related to our series expansion in Section
\ref{sec:compissuesimplissues}. In Section \ref{seccompissues},
details of the method to maximize the marginal likelihood are
given, and in addition we provide computational efficiency gains
and guidelines as to the number of calculations that will occur
using our method. In particular, we initially obtain closed-form
solutions involving infinite sums. Using combinatorial results on
systems of equations with integer coefficients, we show in Theorem
\ref{thm:effunivargamma} how these sums may be re-grouped to a
lengthy (to be discussed and evaluated via simulation) initial
calculation independent of the parameter values, and then a fast
parameter-specific computation which makes the entire approach
tractable. Thus subsequent computations of the marginal likelihood
at different parameter values (necessary for optimization) is
rapid. Additional details of these combinatorial savings are
provided in Appendix \ref{seccombdiophbounds}. In Section
\ref{secnumsim} we provide some simulations to demonstrate the
efficacy of our approach. In Section \ref{subsec:compMCMC} we
compare our closed form series expansion with previous numerical
techniques used to analyze these types of Bayesian inference
problems. In the case of one observation per household, our series
expansions have a comparable run-time to Monte Carlo Markov Chain
methods (in fact, the series expansions are faster); however, for
multiple observations per household these numerical methods are
typically faster, though our series expansions can still be
implemented in a reasonable amount of time. In this manner, our
approach is an alternative, albeit for many practical problems not
one that is faster, but rather one that can be used to verify
other (e.g. MCMC) methods. Some areas for future research and
limitations of our approach, in particular the extension of our
findings to a more general class of priors (Theorem
\ref{thm:compefflincombuvg}), and a more general class of
covariates, are described in Section \ref{secgenfutureresearch}.
We show that, at the cost of introducing new special functions, we
can handle any one-sided probability distribution for the priors.

In Section \ref{sec:covar} we generalize the results of Section
\ref{secindhouseexpand} to allow for covariances. In Section
\ref{subsec:generalG} we derive a closed-form series expansion for
an arbitrary multivariate distribution; however, if the
distribution has a good closed-form expression for its moment
generating function, more can be done. We concentrate on the case
where the response coefficients are drawn from a Multivariate
Gamma distribution, which allows us to have correlations among the
response coefficients. The key observation is Theorem
\ref{thm:MGFmultivariate}, where we interpret the resulting
integrals as evaluations of the Moment Generating Function which
exists in closed-form for the Multivariate Gamma distribution.
Thus we may mirror the arguments from Section
\ref{secindhouseexpand} and again obtain a rapidly convergent
series expansion (Theorem \ref{thm:effmultivargamma}), and the
combinatorial results of Section \ref{seccompissues} and Appendix
\ref{seccombdiophbounds} are still applicable.

Section \ref{conclusionsection} contains some concluding remarks.


\ \\
\begin{center}
\ \ \ \ \section{PROBLEM FORMULATION}\label{problem}
\end{center}

As the logit model and its associated likelihood are well
understood, we briefly describe them in Section \ref{notation} and
focus mainly here (in Section \ref{geometric-series}) on the
geometric series expansion of the model. If we assume the
parameters are independent (all zero covariances), then a
tractable model is obtained by assuming that each is drawn from a
Gamma distribution. This is described in detail in Section
\ref{secindhouseexpand}; in Section \ref{sec:covar} we generalize
the model by assuming the parameters are drawn from a multivariate
Gamma distribution, which allows us to handle covariances among
the parameters. In both cases we obtain closed-form series
expansions. Further, a careful analysis of the resulting
combinatorics leads to computational gains that make the approach
feasible and attractive.

\bigskip
\subsection{Notation}\label{notation}

To describe the model, and to be specific about the data
structures addressed (and not addressed) in this research, we
utilize the following notation. The jargon is drawn from the
Marketing domain and is done for purely explicative purposes.  As
we demonstrate, our approach is applicable for a wide class of
general data structures.

Consider a data set obtained from $i\in \{1,\dots, I\}$ households
(units) containing  $j \in \{1,\dots,J\}$ product categories
(objects; e.g. coffee) measured on $t \in \{1,\dots,N_i\}$
purchase occasions (repeated measures). At each purchase occasion,
for each category $j$ each household $i$ decides whether or not to
purchase in that category.

As is standard, we define
\begin{equation} \twocase{y_{ijt} \ = \ }{1}{if household $i$ buys in category $j$
at time $t$}{0}{otherwise,} \end{equation} \noindent where
$p_{ijt}={\rm Prob}(y_{ijt}=1)$ is the probability of purchase of
the $j$\textsuperscript{th} category by the
$i$\textsuperscript{th} household on its $t$\textsuperscript{th}
purchase occasion. Further, let $P$ denote a set of attributes
describing the categories, with corresponding values $x_{ijt,p}$
such that $X_{ijt}^{T} = (x_{ijt,1}, \dots, x_{ijt,P})$. To
account for differences in base-level preferences for categories,
we define $x_{ijt,1}=1$ defining category-level intercepts.  Thus,
multiplying over all households, categories, and occasions, we
obtain that the standard logit likelihood of the data,
$Y=(y_{ijt})$, is given by
\begin{equation}\label{eqexpandgeo} P(Y|\beta) \ = \ \prod_{i=1}^I
\prod_{j=1}^J \prod_{t=1}^{N_i} \frac{ e^{-X_{ijt}^T \beta_i
y_{ijt} }}{1 + e^{-X_{ijt}^T \beta_i }}, \end{equation} where
$\beta_{i} = (\beta_{i,1},\dots,\beta_{i,P})$ is the coefficient
vector for the $i$\textsuperscript{th} household with variable $p$
specific coefficient, $\beta_{i,p}$.  It is the heterogeneity
across households $i$ in their $\beta_{i,p}$ that we model in
Section \ref{secindhouseexpand} as coming from the Gamma family of
distributions, and in Section \ref{sec:covar} as coming from a
multivariate family of distributions.

The marginalization of the likelihood, which is the problem we
address here, is that we want to ``hit'' $P(Y|\beta)$ (integrate
with respect to) a set of distributions $g(\beta_{i,p}|\Omega)$
depending on parameters $\Omega$ such that
\begin{eqnarray} \label{eq-integrate}
P(Y|\Omega) &  = &  \int P(Y|\beta) g(\beta|\Omega) d \beta
\end{eqnarray}
is available in closed-form.  To accomplish this, we require a
properly chosen series expansion of $P(Y|\beta)$, and we describe
its basic building block next, an application of the geometric
series expansion. The goal is to obtain a good closed-form
expansion of the marginalization of the likelihood for each choice
of parameters $\Omega$. To do so requires finding an appropriate
conjugate distribution leading to tractable integration; we shall
see that the Gamma (Theorem \ref{thm:univargamma}) and
Multivariate Gamma (Theorem \ref{thm:MGFmultivariate})
distributions lead to integrals which can be evaluated in
closed-form. Using such expansions, we then determine the value of
$\Omega$ which maximizes this likelihood; this will allow us to
make inferences about the population heterogeneity distributions.

\subsection{Geometric Series Expansion}\label{geometric-series}

To obtain closed-form Bayesian inferences for the logit model, we
expand the likelihood $P(Y|\beta)$ given in \eqref{eqexpandgeo} by
using the geometric series expansion:
\begin{eqnarray} \label{eq-geometric}
\frac{1}{1 - z} & =  &   \sum_{k=0}^\infty z^k.
\end{eqnarray}
This is directly applicable for our problem as $P(Y|\beta)$ is the
product of terms of the form $\frac{e^{u y}}{1 + e^u}$. Our
interest will be in expanding the denominator when $u < 0$. Note
terms with $u > 0$ can be handled by writing $\frac{1}{1+e^u}$ as
$\frac{e^{-u}}{1+e^{-u}}$.

While in theory we can unify the two cases (positive and negative
values of $u$) by using the $\text{sgn}$ function
($\text{sgn}(u)=1$ if $u>0$, $0$ if $u=0$ and $-1$ otherwise), the
$\text{sgn}$ function is only practical in the case of one
attribute (i.e. $x$ is one dimensional and $P$=1); otherwise it is
undesirable (untenable) in the expansions. In high dimensions
(lots of households with lots of categories and attributes), the
$\text{sgn}$ function leads to numerous, complicated subdivisions
of the integration space. This greatly increases the difficulty in
performing the integration and obtaining tractable closed-form
expansions, and hence is not entertained here. Details of the
unification are available from the authors upon request, and
applying it in practice is an area for future research.  Instead,
we describe below a specific set of restrictions that we employ,
and the class of problems (data sets) where our expansions can then directly
be applied. In Section 6, areas for future research to generalize
our work to richer data sets are discussed.

As mentioned above, to eliminate the need for the sgn function and
to allow for straightforward expansions, we limit our
investigations to the common set of Marketing problems (as
described in Section \ref{notation} and throughout) in which all
\ben \item $X_{ijt,p} \geq 0$, \item $\beta_{i,p} > 0$, \item
$\sum_{p=1}^{P} \beta_{i,p}X_{ijt,p}>0$. \een From a practical
perspective, these restrictions indicate as follows.  First, each
$X_{ijt,p} \geq 0$ is not particularly restrictive, as commonly
utilized descriptor variables -- prices, dummy variables for
feature and display, etc..., as in standard SCANPRO models
(Wittink et al. 1988), are all non-negative and are
straightforwardly accounted for in our framework.  Those variables
which take signs counter to previously signed variables can be
coded as $f(X_{ijt,p})$, for instance $-X_{ijt,p}$ or
$\exp(-X_{ijt,p})$.

Secondly, restriction of $\beta_{i,p} > 0$ may or may not be
restrictive.  If the variables which comprise $X_{ijt,p}$ are ones
in which we want to enforce monotonicity constraints (Allenby,
Aurora, and Ginter, 1995), or naturally one would expect upward
sloping demand at the category-level (which may be much more
likely than at the brand level), then this constraint is not at
all restrictive, and in fact may improve the predictions of the
model. To implement this, especially in the case of dummy coded
$X_{ijt,p}$, the least preferred level should be coded as 0 so
that all other corresponding dummy variables have $X_{ijt,p}=1$
and it is expected that $\beta_{i,p}>0$.

Our third constraint, which is not restrictive as long one of the
$X$'s (e.g. price, coupon, etc...) is non-zero, is required so
that we are not expanding $\frac{1}{2} = \frac{1}{1+1}$ in a
polynomial series; if this condition fails, trivial book-keeping
suffices.

There are two important things to note. First, if all the
$\beta_{i,p}$'s are negative, we may explicitly factor out the
negative sign of each $\beta_{i,p}$, yielding terms such as $-
X_{ijt,p} |\beta_{i,p}|$. If this is the case, for simplicity we
change variables and let $\beta_{i,p} = |\beta_{i,p}|$; thus,
$\beta_{i,p}$ and $X_{ijt,p}$ are now both non-negative, and we
have minus signs in the exponents above. Thus we do not need to
assume all persons have positive or all persons have negative
coefficients, but rather (by recoding $X$ to $-X$) that each
person's coefficients are all positive or all negative. Secondly,
in totality, the restrictions above suggest that our model works
only for categories in which the probability of purchasing in that
category on any given occasion is strictly greater than $\foh$ (if
coded as before) or less than $\foh$ (if coded as in this
paragraph), where again this can very person-by-person. Certainly
an area for future research would be the ability, if possible, to
relax some of these assumptions, and to empirically investigate
the set of product categories for which these restrictions are not
particularly binding (such as long-lasting durable goods).

Thus (after possibly recoding), due to our restrictions, we only
need to use \eqref{eq-geometric} when $X_{ijt}^{T}\beta_{i}
> 0$,  which yields \begin{equation}\frac{1}{1+e^{-X_{ijt}^{T}\beta_{i}}}
 \ \ = \  \ \sum_{k_{ijt}=0}^\infty (-1)^{k_{ijt}}\ e^{-k_{ijt}
X_{ijt}^T \beta_i}, \ \ \ \ X_{ijt}^{T}\beta_{i} \
>\ 0. \end{equation}
This is combined, as described next, with the Gamma or
Multivariate Gamma family of distributions in a conjugate way. It
is the constancy of the sign of $X_{ijt}^{T}\beta_{i}$ that allows
us to use the same series expansion for all $\beta_i$.

\subsection{Expansion of $P(Y|\beta)$}\label{secindhouseexp}

Using the likelihood for the logit model given in (2), we have as
follows:
\begin{eqnarray}\label{eqpybetaexp}  P(Y|\beta)   &  = &  \prod_{i=1}^I \prod_{j=1}^J
\prod_{t=1}^{N_i} \frac{ e^{-X_{ijt}^T \beta_i y_{ijt} }}{1 + e^{-
X_{ijt}^T \beta_i } }  \nonumber \\
 & = & \prod_{i=1}^I \prod_{j=1}^J \prod_{t=1}^{N_i} e^{-y_{ijt}
X_{ijt}^T \beta_i} \prod_{i=1}^I \prod_{j=1}^J \prod_{t=1}^{N_i}
\frac{1}{1 + e^{-X_{ijt}^T \beta_i } } \nonumber \\
 & = & \prod_{i=1}^I \prod_{p=1}^P e^{ - \left(\sum_{j=1}^J
\sum_{t=1}^{N_i} y_{ijt} x_{ijt,p} \right) \beta_{i,p} } \cdot
\prod_{i=1}^I \prod_{j=1}^J \prod_{t=1}^{N_i}
\frac{1}{1 + e^{-X_{ijt}^T \beta_i } } \nonumber \\
 & = & P_1(Y|\beta) P_2(Y|\beta).
\end{eqnarray}
Note that the first term, $P_1(Y|\beta)$, is already an
exponential function. This combines nicely with the Gamma and
Multivariate Gamma distributions, and for each variable
$\beta_{i,p}$, we simply have the exponential of a multiple of
$\beta_{i,p}$. In fact, as we show later in Theorem
\ref{thm:MGFmultivariate}, it is this exponential form that leads
to the result that all closed-form integrals are obtainable using
the Moment Generating Functions of the heterogeneity distribution.
It is the second term, $P_2(Y|\beta)$, that we expand by using the
geometric series. We describe this now.

The real difficulty in coming up with a conjugate family to the
logit model is in the expansion of $P_2(Y|\beta)$. Using the
geometric series expansion, we obtain
\begin{eqnarray} P_2(Y|\beta) & \ = \ & \prod_{i=1}^I \prod_{j=1}^J
\prod_{t=1}^{N_i} \frac{1}{1 + e^{- X_{ijt}^T \beta_i } }
\nonumber\\ & = & \prod_{i=1}^I \prod_{j=1}^J \prod_{t=1}^{N_i}
\sum_{k_{ijt}=0}^{\infty} (-1)^{k_{ijt}} e^{-k_{ijt} X_{ijt}^T
\beta_i } \nonumber\\ & = & \prod_{i=1}^I \prod_{j=1}^J
\prod_{t=1}^{N_i} \sum_{k_{ijt}=0}^{\infty} (-1)^{k_{ijt}}
e^{-k_{ijt} \sum_{p=1}^P x_{ijt,p} \beta_{i,p}}.
\end{eqnarray}
For fixed household $i$, replacing $\prod_{j=1}^J
\prod_{t=1}^{N_i} \sum_{k_{ijt}=0}^\infty$ with
$\sum_{k_{i11}=0}^\infty \cdots \sum_{k_{iJN_i} = 0}^\infty$
yields
\begin{equation}\label{eqpthreetempexpand} P_2(Y|\beta) \ = \  \prod_{i=1}^I
\left( \sum_{k_{i11}=0}^\infty \cdots \sum_{k_{iJN_i}=0}^\infty
(-1)^{\sum_{j=1}^J \sum_{t=1}^{N_i} k_{ijt}} \prod_{p=1}^P  e^{
-\left( \sum_{j=1}^J \sum_{t=1}^{N_i} k_{ijt} x_{ijt,p} \right)
\cdot \beta_{i,p}}\right). \end{equation} Our problem therefore
reduces to finding a good expansion for the integral of
$P_{1}(Y|\beta) P_2(Y|\beta)g(\beta|\Omega)$. We assume
$g(\beta|\Omega)$ is given by a product of Gamma or Multivariate
Gamma distributions (or their generalizations as described below),
rich families of probability densities defined for non-negative
inputs $\beta_{i,p}$ which for certain choices of parameters have
good, closed-form expressions for integrals against exponentials.
For other parameter distribution choices and other probability
distributions, at the cost of introducing new special functions we
still have closed-form series expansions for the integrals (as we
discuss in Section \ref{secgenfutureresearch}). The reason the
Gamma and Multivariate Gamma distributions lead to closed-form
expansions is that both have a good closed-form expression for
their moment generating function; in Theorem
\ref{thm:MGFmultivariate} we generalize our results to any
multivariate distribution with a good closed-form moment
generating function.


\ \\
\begin{center}
\section{UNIVARIATE CASE: GENERALIZED GAMMA}\label{secindhouseexpand} \end{center}

In this section we combine all of the pieces in the case when the
$\beta_{i,p}$ are independently drawn from Gamma distributions
with parameters $(b_p,n_p)$ (independent of $i$): the logit
likelihood given in \eqref{eqexpandgeo}, the geometric series
expansion in \eqref{eq-geometric}, and an integration lemma given
in \eqref{integrationlemma} below which allows us to obtain series
expansions for $P(Y|\Omega)$. Then in Section \ref{seccompissues}
we discuss how to re-group the resulting series expansions for
computational savings. In Section \ref{sec:covar} we consider the
more general case of choosing the $\beta_{i,p}$'s from a
Multivariate Gamma distribution, where now for a given $i$ there
may be correlations among the $\beta_{i,p}$'s.

\subsection{The Gamma Distribution and its
Generalization}\label{sec:gammadistranditsgeneralizations}

As $\beta_{i,p}$ is assumed greater than 0, one flexible
distribution to draw the $\beta_{i,p}$'s from is the three
parameter Generalized Gamma distribution (McDonald and Butler,
1990). The Generalized Gamma distribution is extremely rich, and
by appropriate choices of its parameters, many standard functions
are obtainable. It is defined for $z$ non-negative by
\begin{equation} GG(z;a,b,n) = \frac{|a|}{b \Gamma(n)}
\left(\frac{z}{b}\right)^{an-1} e^{- \left( \frac{z}{b} \right)^a
}. \end{equation} For example, the following assignments of the
parameters $a, b$ and $n$ yield well known distributions:
  $\lim_{a \to 0} GG(z;a,b,n)$ is lognormal;
$GG(z;a,b,1)$ is Weibull;   $GG(z;1,b,n)$ is Gamma; $GG(z;1,b,1)$
is Exponential; $GG(z;2,b,1)$ is Rayleigh. For fixed $a$ and $n$,
the effect of $b$ is to re-scale the units of $z$. That is, as $z$
only appears as $\frac{z}{b}$, $b$ may be interpreted as fixing
the scale (i.e. the commonly interpreted scale parameter); $a$ and
$n$ change the general shape of the Generalized Gamma. Hence, as
opposed to the more familiar Gamma family of distributions
commonly used in Marketing problems, which we focus on here, the
Generalized Gamma has a second shape parameter, $a$, allowing for
more flexible shapes.

We provide in Appendix \ref{subsec:plotsgammas} some plots of the
Gamma family of distributions for various parameter values, and of
a mixture of Gamma distributions, an even more flexible class to
demonstrate its flexibility in providing a rich yet parsimoniously
parameterized set of priors for the $\beta_{i,p}$. Although the
results reported directly in this paper correspond to the
heterogeneity distribution following a single Gamma distribution,
they straightforwardly extend to a mixture of Gammas, where the
mixture is a weighted sum of component Gammas. For each component
of the mixture we can integrate its expansion term by term, and
hence the entire weighted sum. From a practical point of view,
this allows us to handle the situation of latent class modelling,
in which the $\beta_{i,p}$ come from a latent segment, each of
which has its own Gamma parameters.

As the geometric series expansion of the logit likelihood, as
described in Sections \ref{notation} and \ref{geometric-series},
will lead to terms involving exponential functions, our key
integration result arises from integrating an exponential function
against a Gamma distribution. For simplicity we consider only the
case of a Gamma distribution ($a = 1$), and discuss its
generalization below and in Section \ref{secgenfutureresearch}.

This assumption allows us to not only obtain closed-form
expansions, but these expansions will be rational functions of the
arguments of the Gamma distribution (which allow us to obtain
tractable closed-form expansions for the \emph{derivatives} as
well). For notational convenience let $G(z;b,n) = GG(z;1,b,n)$
denote a Gamma distribution with parameters $b$ and $n$.

\begin{lem}[Exponential against a Gamma distribution]\label{lemexpagainstggaisone}

Consider a Gamma distribution $G(z;b,n)$. For $d \ge 0$,
\begin{eqnarray}
 e^{-zd} G(z;b,n)  & = &  (1+bd)^{-n} G\left(z,
\frac{b}{1+bd},n \right).
\end{eqnarray}
As $G(z;\frac{b}{1+bd},n)$ is a probability distribution, we
obtain
\begin{eqnarray} \label{integrationlemma}
 \int_{z=0}^\infty e^{-zd} G(z;b,n) dz  & = & \frac{1}{(1+bd)^n}.
\end{eqnarray}
\end{lem}

\noindent See Appendix \ref{subsectionggintproof} for a proof.
This closed-form integration result allows us to avoid resorting
to Monte Carlo or other numerical techniques to approximate the
integral for $P(Y|\Omega)$, and is our integration ``engine''.
Later in Theorem \ref{thm:MGFmultivariate} of Section
\ref{sec:covar} we generalize Lemma \ref{lemexpagainstggaisone} by
interpreting it as evaluating the moment generating function of
the Gamma distribution at $-d$.

\subsection{Series Expansion for $P(Y|\Omega)$ for the Gamma
Distribution}\label{subsec:seriesexppygamma}

We assume the response coefficients are drawn from the Gamma
distribution. Summarizing our inference problem, we need to
investigate the integral for $P(Y|\Omega)$:
\begin{equation}\label{eqXX}\int_0^\infty \cdots \int_0^\infty P(Y|\beta) g(\beta|\Omega)
d\gb \ = \ \prod_{i=1}^I \int_0^\infty \cdots \int_0^\infty
P_{1i}(Y|\beta) P_{2i}(Y|\beta) \prod_{p=1}^P
G(\beta_{i,p};b_p,n_p) d\beta_{i,p}, \end{equation} where
\begin{eqnarray} P_{1i}(Y|\beta) & \ = \ & \prod_{p=1}^P e^{ -
\left(\sum_{j=1}^J \sum_{t=1}^{N_i} y_{ijt} x_{ijt,p} \right)
\beta_{i,p} } \nonumber\\ P_{2i}(Y|\beta) & \ = \ &
\sum_{k_{i11}=0}^\infty \cdots \sum_{k_{iJN_i}=0}^\infty
(-1)^{\sum_{j=1}^J \sum_{t=1}^{N_i} k_{ijt}} \prod_{p=1}^P  e^{
-\left( \sum_{j=1}^J \sum_{t=1}^{N_i} k_{ijt} x_{ijt,p} \right)
\cdot \beta_{i,p}} \nonumber\\ g(\beta|\Omega) & = & \prod_{p=1}^P
G(\beta_{i,p};b_p,n_p).  \end{eqnarray} Because of the conditional
independence across $i$, we can evaluate each integral in
\eqref{eqXX} separately. We denote each of the $i$-integrals above
by $\inti$ (for the $i$\textsuperscript{th} household), where
\begin{eqnarray}\label{eqinitexpand} \inti & \ = \ & \int_0^\infty \cdots \int_0^\infty
\prod_{p=1}^P \sum_{k_{i11}=0}^\infty \cdots
\sum_{k_{iJN_i}=0}^\infty (-1)^{\sum_{j=1}^J \sum_{t=1}^{N_i}
k_{ijt}} \nonumber\\ & & \ \ \ \cdot \ e^{ - \left(\sum_{j=1}^J
\sum_{t=1}^{N_i} y_{ijt} x_{ijt,p} \right) \beta_{i,p} } e^{
-\left( \sum_{j=1}^J \sum_{t=1}^{N_i} k_{ijt} x_{ijt,p} \right)
 \beta_{i,p}} G(\beta_{i,p};b_p,n_p) d\beta_{i,p} \nonumber\\
 & \ = \ &  \sum_{k_{i11}=0}^\infty
\cdots \sum_{k_{iJN_i}=0}^\infty (-1)^{\oa{k}\cdot\oa{1}}
\prod_{p=1}^P \int_{\beta_{i,p}=0}^\infty e^{ -K_{i,p} \beta_{i,p}
} G(\beta_{i,p};b_p,n_p) d\beta_{i,p},
\end{eqnarray} where \begin{eqnarray}
 K_{i,p} \ = \ \sum_{j=1}^J
\sum_{t=1}^{N_i} (y_{ijt} + k_{ijt})x_{ijt,p}, \ \ \
\oa{k}\cdot\oa{1} \ = \ \sum_{j=1}^J \sum_{t=1}^{N_i} k_{ijt}.
\end{eqnarray} Therefore
\begin{equation}\label{eqXXX} P(Y|\Omega) \ = \ \int_0^\infty \cdots \int_0^\infty P(Y|\beta) g(\beta|\Omega)
d\gb \ = \ \prod_{i=1}^I \inti. \end{equation} Applying the
integration lemma (Lemma \ref{lemexpagainstggaisone}) to
\eqref{eqinitexpand} yields
\begin{eqnarray}\label{eqexpagainstGG}
 \int_{\beta_{i,p}=0}^\infty e^{-K_{i,p}
\beta_{i,p}} G(\beta_{i,p};b_p,n_p)d\beta_{i,p}  & \ = \ &
\frac{1}{\left(1+b_pK_{i,p}\right)^{n_p}}. \end{eqnarray} By
combining the expansion in \eqref{eqinitexpand} with \eqref{eqXX},
we obtain our final result for $\inti$:

\begin{thm}\label{thm:univargamma} Assume the $\beta_{i,p}$ are
independently drawn from Gamma distributions with parameters
$(b_p,n_p)$. Then $P(Y|\Omega) = \prod_{i=1}^I \inti$, where
 \begin{equation}\label{eqinti}  \inti \ =\
\sum_{k_{i11}=0}^\infty \cdots \sum_{k_{iJN_i}=0}^\infty
(-1)^{\sum_{j=1}^J \sum_{t=1}^{N_i} k_{ijt}} \prod_{p=1}^P
\frac{1}{\left(1+b_pK_{i,p}\right)^{n_p}}, \ \ \  K_{i,p} \ = \
\sum_{j=1}^J \sum_{t=1}^{N_i} (y_{ijt} + k_{ijt})x_{ijt,p}.
\end{equation}
\end{thm}

\noindent Hence the log marginal distribution, $\log L = \log
P(Y|\Omega)= \sum_{i}\log(\inti)$, can be computed as the sum of
the logarithm of \eqref{eqinti}. This yields the desired
closed-form solution.

\subsection{Computational and Implementation
Issues}\label{sec:compissuesimplissues}

\subsubsection{Computational Issues and Gains from Diophantine
Analysis}\label{seccompissues}

While Theorem \ref{thm:univargamma} yields a closed-form expansion
for the marginal posterior distribution when the response
coefficients are independently drawn from Gamma distributions, to
be useful we must be able to efficiently determine the optimal
values of the parameters $\Omega$. As written, the number of terms
needed in the series expansions are computationally
expensive/impossible (i.e. the upper bounds are at $\infty$). If
every sum ranged from $0$ to $R$, to have good expansions $R$
would have to be prohibitively large. In this section we describe
a more efficient way to group the summands to significantly
decrease computational time and maximize the marginal posterior
which will make this more computationally tractable for the
Marketing scientist.  We also note that due to the high degree of
non-linearity and the infinite series expansion, there does not
exist a closed-form solution for the optimal parameter values,
$\hat{\Omega}$, by simply solving the first-order condition
equation $\frac{\partial \log L}{\partial \Omega}=0$. We therefore
use numerical methods to obtain the maximum marginal a posteriori
values.

One common approach to determining the optimal values is to use a
multivariate Newton's method. Unfortunately, in many of the
simulations investigated here, the flatness of the surface around
the mode and the multi-modality of the marginal posterior led to
poor convergence; however, we expect for larger (and different)
data sets, Newton's method may become feasible and hence we
include the closed-form first, second, and cross derivatives in
Appendix \ref{secnewtonprogram}.  We also note that one reason for
our choice of the Gamma distribution was that the resulting
expansions (see Theorem \ref{thm:univargamma}) are elementary
functions of the parameters $b_p$ and $n_p$, and hence have
elementary closed-form expansions for their derivatives. This
facilitates calculations of elasticities, shown to be crucial in
determining optimal marketing strategies (Russell and Bolton,
1988).

We therefore instead resorted to evaluating \eqref{eqinti} in a
grid over the parameter space, and then choosing the value that
maximized the marginal posterior.  That is, the beauty and value
of our expansions is that it allows us to calculate $\inti$
rapidly even for many grid points.  However, this is assuming that
we can truncate each of the summations at a computationally
feasible value, an approach we  now describe.

As the expansion stands in \eqref{eqinti}, without careful
thought, only moderate sizes for $J$ and $N_i$ are feasible. In
any numerical calculation, the infinite sums must be truncated.
For simplicity and for explicative purposes, assume each sum
ranges from $0$ to $R-1$. As there are $JN_i$ summations, we have
a total of $R^{JN_i}$ terms to evaluate. Additionally, we have a
product over $p \in \{1,\dots,P\}$ attributes, and then a product
or sum over $i \in \{1,\dots,I\}$ households. If we ssume all $N_i
= N$ for purposes of approximating the computational complexity,
the number of computations required is therefore of order $P
\prod_i R^{JN_i} = PI \cdot R^{JN}$.

As we want to determine the values for the parameters $b_p, n_p$
that maximize the integral given in (18), two common approaches,
Newton's Method or evaluating in a grid, can theoretically be done
(especially as we have explicit formulas); however, the number of
terms makes {\rm direct} computation from this expansion (i.e.
without computational savings as described below) impractical at
present computing speeds. For each parameter, we need to calculate
on the order of $PI\cdot R^{JN}$ terms for just \emph{one}
iteration of Newton's Method or evaluation of $\log L$ for a grid
approach.   We discuss a way to re-group the terms in the
expansion which greatly reduces the computational time and allows
us to handle larger triples $(R,J,N_i)$.

We show below in detail that what allows us to succeed is that it
is possible to re-group the computations in such a way that we
have a lengthy initial computation, whose results we store in a
data file. From this, it is possible to evaluate the
log-likelihood (or derivatives if using Newton's method) at all
points of interest extremely rapidly. The reason such a savings as
described below is possible, in some sense, is that the
computations factor into two components, and most of the
computations are the same for all values of the parameters and
hence only need to be done once.

Consider the case of $P$ attributes: $p \in \{1,2,\dots,P\}$. We
then have $x$-vectors
\begin{equation} x_{i,p} \ = \ (x_{i11,p}, \dots, x_{iJN_i,p}), \ \ p \in
\{1,\dots,P\}. \end{equation} Assume all $x_{ijt,p}$ are integers;
this is not a terribly restrictive assumption\footnote{This is not
restrictive even for an attribute like price. Many studies are
done with a discrete set of integer prices and in other cases,
even if there were a fairly moderate number, the model can handle
it albeit with increased computation.}, and can be simply
accomplished by changing the scale we use to measure the
$x_{ijt,p}$'s. The advantage of having integer $X$'s is that we
now have Diophantine equations, and powerful techniques are
available to count the number of solutions to such equations and
hence ``judge'' the feasible values of $(R,J,N_{i})$.

For notational convenience let $\oa{k} =
(k_{i11},\dots,k_{iJN_i})$, $\oa{1} = (1,\dots,1)$ and
\begin{eqnarray}\label{eqykip} Y_{i,p} & \ = \ & \sum_j\sum_t y_{ijt}x_{ijt,p}
\nonumber\\ K_{i,p} & \ = \ & \sum_{j=1}^J \sum_{t=1}^{N_i}
(y_{ijt} + k_{ijt})x_{ijt,p} \ = \ Y_{i,p} + \sum_{j=1}^J
\sum_{t=1}^{N_i} k_{ijt}x_{ijt,p}.\end{eqnarray} Recall from
\eqref{eqinti} that when the $\beta_{i,p}$ are independently drawn
from Gamma distributions that
 \begin{equation}\label{eqnewinti} \inti \ =
\sum_{k_{i11}=0}^\infty \cdots \sum_{k_{iJN_i}=0}^\infty
(-1)^{\sum_{j=1}^J \sum_{t=1}^{N_i} k_{ijt}} \prod_{p=1}^P
\frac{1}{\left(1+b_pK_{i,p}\right)^{n_p}}.\end{equation} Fix an $i
\in \{1,\dots,I\}$. We see that \eqref{eqnewinti} depends weakly
on $\oa{k}$; by \eqref{eqykip}, all that matters are the dot
products $\oa{k}\cdot \oa{x}_{i,p}$ and the parity of
$(-1)^{\oa{k}\cdot\oa{1}}$. Let $r = (r_1,\dots,r_P)$ be a
$P$-tuple of non-negative integers. For each $\oa{k}$ we count the
number of solutions to the system of Diophantine equations
$\oa{k}\cdot\oa{x}_{i,p} = r_p$ ($p \in \{1,\dots,P\}$) while
recording the sign of $(-1)^{\sum_j\sum_t k_{ijt}} =
(-1)^{\oa{k}\cdot\oa{1}}$. Explicitly, we may re-write
\eqref{eqinti} from Theorem \ref{thm:univargamma} as

\begin{thm}\label{thm:effunivargamma}
Set $Y_{i,p}  =  \sum_j\sum_t y_{ijt}x_{ijt,p}$ and
\begin{eqnarray}\label{eq:dpcount} S(M) & \ = \ & \left\{v: v =
(v_1,\dots,v_{M}), v_l \in \{0, 1, 2, 3, \dots\} \right\}
\nonumber\\  K_i(x,r,+) & \ = \ & \#\{k \in S(JN_i): \forall p \in
\{1,\dots,P\}, \oa{k}\cdot\oa{x}_{i,p} = r_p,
(-1)^{\oa{k}\cdot\oa{1}} = +1 \} \nonumber\\ K_i(x,r,-) & \ = \ &
\#\{k \in S(JN_i): \forall p \in \{1,\dots,P\},
\oa{k}\cdot\oa{x}_{i,p} = r_p, (-1)^{\oa{k}\cdot\oa{1}} = -1 \}.
\end{eqnarray}
Assume the $\beta_{i,p}$ are independently drawn from Gamma
distributions with parameters $(b_p,n_p)$. Then $P(Y|\Omega) =
\prod_{i=1}^I \inti$ with
\begin{equation}\label{eqintigoodone}\inti \ = \ \sum_{r \in S(P)}
\frac{K_i(x,r,+)-K_i(x,r,-)}{\prod_p (1 + b_pY_{i,p} +
b_pr_p)^{n_p}}.
\end{equation}
\end{thm}

There is a large computational startup cost in solving
\eqref{eq:dpcount}, but future computations are significantly
faster. We calculate $K_i(x,r,\pm)$ \emph{once}, and store the
results in a data file. Then, in subsequent calculations, we need
only input the new values for $b_p$ and $n_p$ (or even better
calculate the values for multiple $b_p$ and $n_p$ simultaneously
if we are evaluating over a grid). The advantage of such an
expansion is that successive terms involving larger $r$ decay with
$\oa{k}\cdot\oa{x}_{i,p}$. Thus, we do not want to truncate the
sum $\sum_{k_{ij1}} \cdots \sum_{k_{ijP}}$ by having each sum
range from $0$ to $R-1$; instead, we want to consider the
$k$-tuples where the dot products are small, as the $k$-dependence
is weak (the expansion depends only on the value of
$\oa{k}\cdot\oa{x}_{i,p}$). From a computational point of view,
there is enormous savings in such grouping.

To determine how many terms are needed for this truncation to be a
good approximation to the infinite expansion requires an analysis
of $K_i(x,r,+) - K_i(x,r,-)$. We sketch some straightforward,
general bounds in Appendix \ref{seccombdiophbounds}. We do not
exploit the gain from the factors of $(1 + b_pY_{i,p} +
b_pr_p)^{-n_p}$ so that our bounds will apply to the more general
cases that we consider later (explicitly, the multivariate
distributions with good closed-form moment generating functions of
Section \ref{sec:covar}). One other point to note and which will
greatly improve the convergence of the truncated expansions is to
introduce translations in the Gamma distributions, which will give
exponentially convergent factors. Assume each $\beta_{i,p} \ge
\epsilon$; for $\gep$ small (e.g. 0.0001), from a practical point
of view such an assumption is harmless as a coefficient restricted
to this range is not practically different than one restricted to
be greater than or equal to 0. Explicitly, we draw $\beta_{i,p}$
from $G(z-\epsilon;b_p,n_p)$ rather than $G(z;b_p,n_p)$. Similar
arguments as before yield
\begin{equation}\label{eqgoodYPexpansion} \inti \ = \ \sum_{r
\in S(P)} \prod_{p=1}^P e^{-Y_{i,p} \gep}  \frac{
(K_i(x,r,+)-K_i(x,r,-)) e^{-r_p\gep}}{ (1 + b_pY_{i,p} +
b_pr_p)^{n_p}}.
\end{equation} As $K_i(x,r,+)-K_i(x,r,-)$ grows at most
polynomially (see Theorem \ref{thm:polybound} in Appendix B), it
is clear the above expansion converges (and for reasonable values,
it will converge more rapidly than when $\epsilon=0$).

\subsubsection{Numerical Simulations}\label{secnumsim}

To demonstrate the efficacy of our approach given in
\eqref{eqnewinti}, we ran a series of numerical simulations.  The
results reported here are from two sets of the many simulations
conducted, the remainder of which are available upon request. The
first simulation design was chosen to be computational feasible;
however, without loss of generality it contains all the elements
that are required to generalize our results.  In some sense, due
to its maximal sparseness in information, it is the most strict
test of our approach.

Specifically, we report here first on a series of simulations with the following design:
\begin{itemize}
\item{$P=1$, one attribute per observation,} \item{$JN_{i}=1$, one
brand and one observation per household,} \item{$I=1000$, one
thousand households, } \item{$0 \leq k_i \leq R$ with $R$ (the
number of polynomial expansion terms) equal to 100,} \item{$a=1$
(the Gamma distribution) for various choices of $b$ and $n$, and}
\item{untranslated Gamma distribution (i.e. $\epsilon = 0$ in
\eqref{eqgoodYPexpansion}).}
\end{itemize}

\noindent All simulations were run using Matlab on a 1.9 GHZ
athlon processor with 192 Mb of RAM, a very modest computing
machine in today's standards.

For each $b_1$ and $n_1$ pair, 25 simulates were run by: (i)
choosing $I=1000$ values of $\beta_{ip=1}$ from a $G(z;b_1,n_1)$.
The values of $x_{ijt}$ were selected from the values $(1,2,3)$
with equal probability, and then arbitrarily scaled by a constant
$c$ to make the values of $\beta_{i,p} \cdot x_{ijt}$ reasonable
so as to allow for enough 0/1 variation in the $y_{ijt}$.  Then
for each of the 25 simulates, we numerically approximated
$P(Y|\Omega)$ as given by \eqref{eqnewinti} and then maximized the
resulting marginal likelihood, as a function of $\Omega$, using a
grid of values.  In particular, we utilized a grid size of
dimension $5\times 7$ centered at $(b_1,n_1)$ with spacings of
$.1$ (this was reached after considerable empirical testing to
ensure enough fineness and that the solutions were not occurring
on the boundary of the grid).


We summarize our results in the table below: the true values of
$b_1$ and $n_1$, the mean and standard deviation over 25
replicates of the estimated values, and the $t$-statistics for
both $b$ and $n$.

\begin{center}
\begin{tabular}{|r||r|r||r|r|}
  $(b_1,\ n_1)$ & $(\overline{b_1},\ \overline{n_1})$ & $(\sigma_{b_1},\ \sigma_{n_1})$
  & t-stat ($b_1$) & t-stat ($n_1$) \\
    \hline
(5,\ 14)   &   (5.21,\ 14.86)    &   (1.40,\ 3.36) &
0.76    &   1.29 \\

(10,\ 28) & (10.72,\ 26.38) &  (1.46,\ 3.17) &    2.46 &   -2.55 \\

\hline

(9,\ 9)  &   (9.06,\ 9.64) & (2.41,\  2.37)    & 0.12 &
1.36 \\

(18,\   18)  &   (17.39,\ 18.62)   &   (2.28,\ 2.40) & -1.34
&   1.30\\

\hline

(11.5,\  6.5) &  (10.65,\ 7.38) &  (2.46,\ 2.03)   & -1.73 &
2.16 \\

(23,\  13) &   (23.93,\ 12.63)   &   (2.35,\ 1.43) &
1.97    &   -1.30\\

\end{tabular}
\end{center}

To assess whether the simulate values are in accordance with the
true values, we conducted t-tests for each of the parameters and
simulated conditions.  This resulted in 12 significance tests, all
of which correspond to a t-distribution with 24 degrees of freedom
(note we did 25 simulates).  Using the common, albeit
conservative, Bonferroni adjustment method for multiple
comparisons, we note the critical value of 3.167 in absolute value
of which none of the comparisons is close (the corresponding value
for one comparison is 2.064, which 9 of the 12 are less than).
This suggests a very adequate fit of our approach and therefore
the size of $R$ in our polynomial expansions.  Other simulations,
not shown, suggested higher values of $R$ provided even greater
accuracy.

The six set of simulations were chosen to be indicative of three
possible settings, $b_1>n_1, b_1=n_1$ and $b_1<n_1$.  We then
replicated these three settings by scaling each of the values of b
and n by a factor of 2.  In this way we are able to show that it
is not a particular ordering of $b_1$ and $n_1$ that matters nor
the relative sizes of them.  Note again, as above, that this
simulation test of our approach is ultra-conservative in that we
have tested our method using $J \cdot N_i = 1$ and $I=1000$,
modest values.  That is, with simply one observation per household
and 1000 households, our approach is accurately able to
reconstruct the heterogeneity distribution from which the
$\beta_{i,p}$ were derived. This result was also not dependent on
$I$=1000, as shown below, and hence would have led to even faster
processing time. Our belief is that this is a strong signal of the
efficacy of our approach.

A second series of simulations with more general conditions was
conducted in which the number of attributes was increased to
$P=2$. The purpose was to see how well {\em multiple} Gamma
distributions could be detected. There are now four parameters
($b_1,n_1,b_2,n_2$). To have these simulations run in a comparable
time as the previous, we chose $I = 250$, $R = 40$, a grid of size
$4\times 4 \times 4 \times 4$ centered at $(b_1,n_1,b_2,n_2)$ with
a grid spacing of $.5$ units, and 10 simulates for each condition.
We summarize the results below.

\begin{center}
\begin{tabular}{|r||r|r||r|r|r|r|}
  $(b_1,n_1,b_2,n_2)$ & $(\overline{b_1},\overline{n_1},\overline{b_2},
  \overline{n_2})$ &
  $(\sigma_{b_1},\sigma_{n_1},\sigma_{b_2},\sigma_{n_2}$)
  & $\text{t-stat} \atop (b_1)$ & $\text{t-stat} \atop (n_1)$ &
  $\text{t-stat} \atop (b_2)$ & $\text{t-stat} \atop (n_2)$ \\
    \hline
(9,\ 9,\ 18,\ 18)   &   (8.60,\ 8.75,\ 17.9,\ 18.1)  & (2.12,\
2.20,\ 1.96,\ 2.25)
& -.60 & -.36 & -.16 & .14 \\

(11.5,\ 6.5,\ 23,\ 13) & (11.3,\ 6.80,\ 22.3,\ 12.8) & (1.77,\
1.95,\ 2.08,\ 1.96) & -.36 & .49 & -1.14 & -.40\\

(5,\ 14,\ 23,\ 13) & (4.85,\ 13.9,\ 24.15\, 14.05) & (1.56,\
1.76,\ 1.55,\ 1.28) & -.30 & -.18 & 2.35 & 2.60




\end{tabular}
\end{center}

This resulted in 12 significance tests, all which correspond to a
t-distribution with 9 degrees of freedom (note we did 10
simulates).  Using the common, albeit conservative, Bonferroni
adjustment method for multiple comparisons, we note the critical
value of 3.81 in absolute value of which none of the comparisons
is close (the corresponding value for one comparison is 2.26,
which ten of the twelve values are less than). This suggests a
very adequate fit of our approach and therefore the size of $R$ in
our polynomial expansions.

Our findings suggest again the general efficacy of our approach as
none of the significance tests indicate divergence between the
true and estimated parameter values.

\subsubsection{Comparison with Monte Carlo Markov Chain
Methods}\label{subsec:compMCMC}

As mentioned previously, and described in detail in Appendix
\ref{seccombdiophbounds}, one aspect of our theoretical results
that requires study is its computational feasibility due to the
large number of summands. As the {\em exact} results in Theorem
\ref{thm:univargamma} or Theorem \ref{thm:effunivargamma} have
upper sum limits at infinity, we conducted an additional small
scale simulation to assess the efficacy of our method under the
truncation approximation. To act as a further baseline to our
approach, we also ran a Bayesian MCMC sampler (a Bayesian
multinomial logit model with non-conjugate gamma priors as per
Section \ref{secindhouseexpand}) to assess both the computation
accuracy for our approach and its accuracy per unit time compared
to established extant methods. All analyses were run on a Pentium
IV 3.3MHZ processor with 2GB of RAM. For a more accurate
comparison of times we wrote a C program rather than a Matlab
program (as in \ref{secnumsim}) for evaluating the truncated sums.

In particular, we simulated data for $I = 1000$ households, $N_{i}
= 1$ or $5$ observations per household, generated by a multinomial
logit model (see \eqref{eqexpandgeo}) with $P = 1$ covariates.
Each household's value of $\beta_{i}$ was drawn from a Gamma
distribution\footnote{Computation time was essentially invariant
over the exact values of $b$ and $n$ chosen. The run-time is a
polynomial in $R$ of degree $N_i$; further empirical testing is
needed to ascertain how well our approach works in these
settings.} with $b=5$ and $n=14$. To analyze our approach, we
evaluated the approximated marginal likelihood (marginalized over
$\beta_{i}$) over a grid (of size $21 \times 21$) using the
Diophantine computation savings by only looking at sums with
$k_{1}+ \cdots + k_{N_{i}} \le R$ for various choices of $R$ (as
compared to $N_{i}$ sums where each went from 0 to $R$, which
leads to the inclusion of many summands of negligible size).

When $N_i = 5$ the Bayesian MCMC sampler for 6000 iterations for 3
chains (0.01667 seconds per iteration) took about 50 seconds,
where the convergence diagnostic of Gelman and Rubin (1992)
indicated convergence after approximately 3000 iterations (hence
9000 observations available for estimation after burn-in). For
$N_i=1$ the Bayesian MCMC sampler for 6000 iterations for 3 chains
(0.01667 seconds per iteration) took about 20 seconds.

For the C program based on our truncated series expansions, the
approximations depend on the parity of $R$ (if $R$ is even then
the final summands all have a factor of $+1$, while if $R$ is odd
the final summands all have a factor of $-1$). Thus if the
resulting values at the grid points are stable for two consecutive
values of $R$, we have almost surely included enough terms in our
truncation. For $N_i = 1$ there was about a 2\% difference in
values when $R=100$ and $101$ (about 12 seconds); there was about
a .2\% difference in values when $R=200$ and $201$ (about 24
seconds). These run-times compare favorably with those of the
Bayesian MCMC sampler. For more observations per household,
however, the Bayesian MCMC sampler does better. The problem, as
shown in Appendix \ref{seccombdiophbounds}, is that the number of
summands with $k_1 + \cdots + k_{N_i} \le R$ is a polynomial in
$R$ of degree $N_i$. When $N_i = 5$ and $R=6$ the program ran for
about 40 seconds, and when $N_i = 5$ and $R=7$ the run-time was
about 64 seconds; while these values of $R$ are too small to see
convergence in the truncated series, for these data sets the
series expansion is still implementable, though at a cost of a
significantly greater run-time.

Thus our series expansions, with the present computing power, are
comparable to existing numerical methods only in the case of one
observation per household, though they can still be implemented in
a reasonable amount of time for multiple observations.


\subsection{Generalizations of the Univariate Gamma
Distribution}\label{secgenfutureresearch}

We describe several natural generalizations of our model. We
assume for each $i$ that $\beta_{i,1},\dots,\beta_{i,P}$ are
independent below; see Section \ref{sec:covar} for removing this
assumption as well.

At the expense of using special functions, we may easily remove
the assumption that the $\beta_{i,p}$ are drawn from a Gamma
distribution; however, as the research currently stands, we still
must assume the $\beta_{i,p}$ are drawn from one-sided
distributions. In the case of just one attribute, it is
straightforward to generalize our methods to handle $\beta_{i,p}$
drawn from \emph{any} distribution (we split the integration into
three parts, $\beta_{i,p} < \gep, |\beta_{i,p}| \le \gep,
\beta_{i,p}
> \gep$); a natural topic for future research is to handle
$\beta_{i,p}$ drawn from two-sided distributions with multiple
attributes.

\subsubsection{Weakening $\bip > \gep$}

The assumption that $\bip > \gep$ is problematic if we desire to
test the hypothesis that $\bip = 0$. To this end, for each $\bip$
we consider instead of an $\epsilon$-translated Gamma distribution
a 0-point mass Gamma Distribution given by
\begin{equation}w_p \delta(\bip) + (1 - w_p) G(\bip-\gep;b_p,n_p). \end{equation}
In the above, $\delta(x)$ is the Dirac Delta Functional with unit
mass concentrated at the origin; $w_p \in [0,1]$ is a weight and
can be interpreted as a ``weight of evidence'' for
$\beta_{i,p}=0$. It is easier, though by no means necessary, to
obtain closed-form integrals if we assume instead that we have
\begin{equation}w_p \prod_{i=1}^I \delta(\bip)d\bip \ + \ (1-w_p)
\prod_{i=1}^I G(\bip-\gep;b_p,n_p) d\bip. \end{equation} That is,
for a given attribute, either $\bip = 0$ for \emph{all}
households, or they are \emph{all} drawn from a translated Gamma
distribution.

Note, we now have {\em either} translated Gamma distributions
\emph{or} delta masses. If everything were a delta mass, we would
be left with $(-1)^{\oa{k}\cdot\oa{1}}$. In this case, we would
not use the geometric series expansion, as the integration is
trivial.

The expansions are more involved if we have some delta masses and
some non-delta masses (varying across attributes). We would have
to go through the same arguments as above to estimate convergence,
but instead of having $P$ terms in the exponentials, we would have
$P-1$, $P-2$, and so on.

A stronger assumption, leading to the easiest integration, is the
following:
\begin{equation}w \prod_{i=1}^I \prod_{p=1}^P \delta(\bip)d\bip  \ + \ (1-w)
\prod_{i=1}^I \prod_{p=1}^P G(\bip-\gep;b_p,n_p)d\bip.
\end{equation}
That is, \emph{either} everything is from a delta mass, \emph{or}
everything is from some translated Gamma distribution, with a
translation of $\gep$.  In this instance, our approach can be
directly applied.

\subsubsection{Linear Combinations of Gamma
Distributions}\label{subsecgeneralize}

We can increase the flexibility of the model by considering linear
combinations of Gamma distributions:
\begin{equation}\label{eq:betalincombuvg}
w_{p,1} G(\beta_{i,p}-\gep;b_{p,1},n_{p,1}) + \cdots + w_{p,C}
G(\beta_{i,p}-\gep;b_{p,C},n_{p,C}), \end{equation} where
\begin{equation}\forall p: w_{p,1} + \cdots + w_{p,C} \ = \ 1, \ \ \ w_{p,c}
\in [0,1]. \end{equation} We can regard the weights as either new,
additional parameters, or fixed, and \eqref{eqexpagainstGG}
becomes
\begin{equation}\int_{\beta_{i,p}=0}^\infty e^{-K_{i,p} \beta_{i,p}} \sum_{c=1}^C
w_{p,c} G(\beta_{i,p}-\gep;b_{p,c},n_{p,c})d\beta_{i,p} \ = \
\sum_{c=1}^C \frac{w_{p,c}\ e^{-K_{i,p}\gep}}{ \left(1+b_{p,c}
K_{i,p}\right)^{n_{p,c}}}. \end{equation} The essential point is
that, in the above integration, $K_{i,p}$ \emph{does not depend on
$c$}. Thus, we will still have the computational savings, and need
only count the solutions to the Diophantine system once. The
difference is we now have more terms to evaluate, but we still
have rapid savings, and \eqref{eqgoodYPexpansion} becomes

\begin{thm}\label{thm:compefflincombuvg} Notation as in Theorem
\ref{thm:effunivargamma}, let the $\beta_{i,p}$ be independently
drawn from linear combinations of Gamma distributions as in
\eqref{eq:betalincombuvg}. Then $P(Y|\Omega) = \prod_{i=1}^I
\inti$ with
\begin{equation}\inti \ = \ \sum_{c=1}^C \sum_{r \in S(P)}
\prod_{p=1}^P e^{-Y_{i,p} \gep} \frac{w_{p,c}
(K_i(x,r,+)-K_i(x,r,-))e^{-r_p\gep}}{ (1 + b_{p,c}Y_{i,p} +
b_{p,c}r_p)^{n_{p,c}}}.
\end{equation}
\end{thm}

\subsubsection{More General One-Sided Distributions}

There is no a priori reason or necessity to choose $\beta_{i,p}$
from  a  Gamma distribution $G(\beta_{i,p};b_p,n_p)$ (or linear
combinations of these). Because we were assuming the $\beta_{i,p}
\ge 0$ (later, when we wanted $\beta_{i,p} \ge \gep$, this merely
caused us to study translated Gamma distributions), it is natural
to choose a one-sided, flexible distribution such as the Gamma
distribution. If we take any one-sided distribution and translate,
we obtain a similar formula as in  \eqref{eqintigoodone} or
\eqref{eqgoodYPexpansion}. The only difference would be the
functional form of the non-Diophantine piece. The exponential
decay (arising from the requirement that $\beta_{i,p} \ge \gep$)
is still present; it came solely from the geometric series
expansions.

As we have not been using properties of the integration of an
exponential against a Gamma distribution to obtain our convergence
bounds, our arguments are still applicable; however, in general we
don't have simple closed-form expansions with elementary
functions. At the cost of introducing new special functions, we
could handle significantly more general one-sided distributions.
Our integration lemma (Lemma \ref{lemexpagainstggaisone}) is
trivially modified, and we still have computational savings. As we
shall see in Theorem \ref{thm:MGFmultivariate}, our method is
directly applicable to any distribution (univariate or
multivariate) with a closed-form moment generating function.

There are two costs. The first is the introduction of new special
functions in the expansions of the $\inti$; however, by tabulating
these functions once, subsequent evaluations can be done
efficiently. The second difficulty is that, if one attempts to use
Newton's Method, closed-form elementary expansions of the
derivatives are no longer available in many cases; for cases where
the expansions exist, one must calculate the partial derivatives
in a manner similar to that in Appendix \ref{secnewtonprogram}
(for the Gamma distribution).

\ \\
\begin{center}
\ \section{INCORPORATING COVARIANCES: THE MULTIVARIATE GAMMA
MODEL}\label{sec:covar}
\end{center}

Our previous investigations have assumed that the households'
$\beta_{i,1},\dots,\beta_{i,P}$ are independent and that
$\beta_{i,1}, \dots, \beta_{i,P}$ are independently drawn from
Gamma distributions (with different parameters for each
$\gb_{i,p}$). In \S\ref{secgenfutureresearch} we have seen how to
generalize to the case when the $\gb_{i,p}$ are still independent
but drawn from other distributions. We now discuss another
generalization, namely removing the independence assumption of the
$\gb_{i,p}$ and thus allowing non-zero covariances.

Let us assume that the households are still independent; however,
$\gb_{i,1}, \dots, \gb_{i,P}$ are no longer assumed to be
independent. Let us assume that for each household these are drawn
from some distribution \be G(\gb_{i,1},\dots,\gb_{i,P};\oa{b}) \ =
\ G(\beta_i;\oa{b}), \ee where $\oa{b}$ is some set of parameters.
Our previous work in Section \ref{secindhouseexpand} is the case
\bea\label{eq:previous} \oa{b} & \ = \ &
(b_{i,1},\dots,b_{i,P}, n_{i,1}, \dots, n_{i,P}) \nonumber\\
G(\gb_i;\oa{b}) & = & \prod_{i=1}^P \frac1{\Gamma(n_{i,p})}
\left(\frac{\gbip}{b_{i,p}}\right)^{n_{i,p}} \ e^{-\gbip /
b_{i,p}}. \eea By using a multivariate distribution we can capture
correlations between the coefficients of different brands (the
univariate distribution of \eqref{eq:previous} has all covariances
zero), or in general the coefficients of the covariates.

As we no longer assume that $G$ factors into distributions for
each $\gb_{i,p}$, \eqref{eqinitexpand} is no longer valid and we
now must analyze, for each household $i$,
\begin{eqnarray}\label{eqinitexpandGG} \inti & \ = \ & \int_0^\infty \cdots \int_0^\infty
\sum_{k_{i11}=0}^\infty \cdots \sum_{k_{iJN_i}=0}^\infty
(-1)^{\sum_{j=1}^J \sum_{t=1}^{N_i} k_{ijt}} \nonumber\\ & & \ \ \
\cdot \ \prod_{p=1}^P e^{ - \left(\sum_{j=1}^J \sum_{t=1}^{N_i}
(y_{ijt} + k_{ijt}) x_{ijt,p} \right) \beta_{i,p} }
G(\beta_i;\oa{b}) d\beta_{i,p} \nonumber\\
 & \ = \ &  \sum_{k_{i11}=0}^\infty
\cdots \sum_{k_{iJN_i}=0}^\infty (-1)^{\oa{k}\cdot\oa{1}}
\int_{0}^\infty \cdots \int_{0}^\infty e^{ -K_{i,1} \beta_{i,1} }
\cdots e^{ -K_{i,P} \beta_{i,P} } \cdot G(\beta_i;\oa{b})
d\beta_{i,1} \cdots d\beta_{i,P},
\end{eqnarray} where as before \begin{eqnarray}
 K_{i,p} \ = \ \sum_{j=1}^J
\sum_{t=1}^{N_i} (y_{ijt} + k_{ijt})x_{ijt,p}, \ \ \
\oa{k}\cdot\oa{1} \ = \ \sum_{j=1}^J \sum_{t=1}^{N_i} k_{ijt}.
\end{eqnarray}

Of course, for general $G$ it will be difficult to evaluate
\eqref{eqinitexpandGG} in a tractable form for numerical
computation. One of the advantages of our previous method is that
the integral of an exponential and a gamma distribution was
another gamma distribution, and thus the integrals which arose
were simple expressions of the parameters.

There are two natural ways to proceed. For a general multivariate
distribution $G$ we will be unable to develop a closed-form
expression for the integral in \eqref{eqinitexpandGG} that is
analogous to the one we found for the case of the $\gbip$'s
independently drawn from Gamma distributions (Lemma
\ref{lemexpagainstggaisone}). Instead we could series expand the
remaining exponentials, recognizing the resulting integrals as the
moments of the multivariate distribution.

Alternatively, if $G$ has a known closed-form expression for its
moment generating function, then we may recognize
\eqref{eqinitexpandGG} as simply evaluating this moment generating
function at $(t_1,\dots,t_P) = (-K_{i,1},\dots,-K_{i,P})$.
Unfortunately sometimes the moment generating functions only exist
for suitably restricted $(t_1,\dots,t_P)$, in which case we
combine this approach with the series expansion for the remaining
$P$-tuples. We present these details below.

\bigskip
\subsection{Series Expansion for $P(Y|\Omega)$}

\subsubsection{General Multivariate $G$}\label{subsec:generalG}

For each of the $P$ exponential terms $e^{ -K_{i,p} \beta_{i,p} }$
we may expand in a geometric series, \be e^{ -K_{i,p} \beta_{i,p}
} \ = \ \sum_{\ell_p = 0}^\infty \frac{(-K_{i,p}
\beta_{i,p})^{\ell_p}}{\ell_p!}. \ee Thus \eqref{eqinitexpandGG}
becomes
\begin{eqnarray}\label{eqinitexpandGGbb} \inti &\ =\ &
\sum_{k_{i11}=0}^\infty \cdots \sum_{k_{iJN_i}=0}^\infty
(-1)^{\oa{k}\cdot\oa{1}} \int_{0}^\infty \cdots \int_{0}^\infty
e^{ -K_{i,1} \beta_{i,1} } \cdots e^{ -K_{i,P} \beta_{i,P} } \cdot
G(\beta_i;\oa{b}) d\beta_{i,1} \cdots d\beta_{i,P} \nonumber\\ & =
& \sum_{k_{i11}=0}^\infty \cdots \sum_{k_{iJN_i}=0}^\infty
(-1)^{\oa{k}\cdot\oa{1}} \nonumber\\ & & \ \ \ \ \ \ \
\int_{0}^\infty \cdots \int_{0}^\infty
\sum_{\ell_1,\dots,\ell_P=0}^\infty
\frac{(-K_{i,1}\gb_{i,1})^{\ell_1}}{\ell_1!} \cdots
\frac{(-K_{i,P}\gb_{i,P})^{\ell_P}}{\ell_P!}\ G(\beta_i;\oa{b})
d\beta_{i,1} \cdots d\beta_{i,P} \nonumber\\ &=&
\sum_{k_{i11}=0}^\infty \cdots \sum_{k_{iJN_i}=0}^\infty
(-1)^{\oa{k}\cdot\oa{1}} \sum_{\ell_1,\dots,\ell_P=0}^\infty
\frac{(-K_{i,1})^{\ell_1} \cdots (-K_{i,P})^{\ell_P}}{\ell_1!
\cdots \ell_P!} \mu_{\ell_1,\dots,\ell_P},
\end{eqnarray} where \be \mu_{\ell_1,\dots,\ell_P} \ = \
\int_0^\infty \cdots \int_0^\infty \beta_{i,1}^{\ell_1} \cdots
\beta_{i,P}^{\ell_P} \cdot G(\beta_{i,1},\dots,\beta_{i,P};\oa{b})
d\gb_{i,1}\cdots d\gb_{i,P}. \ee

We thus obtain a closed-form expression again, except now we have
additional summations over $\ell_1, \dots, \ell_P$. Here
$\mu_{\ell_1,\dots,\ell_P}$ is the $(\ell_1,\dots,\ell_P)$
non-centered moment of the distribution $G$. For a general
distribution these may be difficult to evaluate explicitly; we
need a one-sided distribution (with some parameters $\oa{b}$) that
is flexible in terms of shape as well as having good formulas for
the moments $\mu_{\ell_1,\dots,\ell_P}$.

Our combinatorial results from Section \ref{seccompissues} (where
we were able to re-arrange calculations to save computational
time) depended crucially on the fact that the exponential versus
gamma integrals from before led to simple expansions such as $(1 +
b_pK_{i,p})^{-n_p}$; these expansions did not depend on the actual
values of $k_{i11},\dots, k_{iJN_i}$ but only some linear
combinations (dot products). Thus we still have combinatorial
savings in the $k_{ijt}$ sums.

\subsubsection{Multivariate $G$ with Closed Form Moment Generating
Functions}\label{subsec:goodmgfG}

Let $\beta_i = (\gbio,\dots,\gbiP)$ be distributed according to a
multivariate density $G(\gbi,\oa{b})$. The moment generating
function of $G$ is given by \be M_{\beta_i}(t_1,\dots,t_P) \ = \
\E\left[ e^{t_1 \beta_{i,1} + \cdots + t_P \beta_{i,P}}\right],
\ee where the expectation is with respect to $G$; i.e., \be
M_{\beta_i}(t_1,\dots,t_P) \ = \ \int_{\gbio} \cdots \int_{\gbiP}
e^{t_1 \beta_{i,1} + \cdots + t_P \beta_{i,P}} \
G(\gbio,\dots,\gbiP;\oa{b}) d\gbio \cdots d\gbiP. \ee Depending on
the distribution, the moment generating function may exist for all
$P$-tuples $(t_1,\dots,t_P)$, or instead only for suitably
restricted $P$-tuples. We immediately obtain

\begin{thm}\label{thm:MGFmultivariate}
Assume the moment generating function $M_{\gbi}(t_1,\dots,t_P)$
for the multivariate distribution $G(\beta_i,\oa{b}) =
G(\beta_{i,1},\dots,\beta_{i,P},\oa{b})$ exists for all
$(t_1,\dots,t_P)$. Then $P(Y|\Omega) = \prod_{i=1}^I \inti$, where
\be H_i \ = \ \sum_{k_{i11}=0}^\infty \cdots
\sum_{k_{iJN_i}=0}^\infty (-1)^{\oa{k}\cdot\oa{1}}
M_{\gbi}(-K_{i,1},\dots,-K_{i,P}), \ee with  \be K_{i,p} \ = \
\sum_{j=1}^J \sum_{t=1}^{N_i} (y_{ijt} + k_{ijt})x_{ijt,p}, \ \ \
\oa{k}\cdot\oa{1} \ = \ \sum_{j=1}^J \sum_{t=1}^{N_i} k_{ijt}. \ee
\end{thm}

\subsection{Multivariate Gamma
Distributions}\label{subsec:specialmultdistr}

We list several versions of Multivariate Gamma distributions (with
non-zero covariances) that have closed-form expressions for their
moment generating functions, and thus satisfy the conditions of
Theorem \ref{thm:MGFmultivariate}. For additional multivariate
distributions see Appendix \ref{sec:moremultdist}. All page and
equation references in Sections \ref{sec:cherrama} and
\ref{subsec:mvgdqazxsw} and are from Kotz, Balakrishnan and
Johnson 2000.

\subsubsection{(Cheriyan and Ramabhadran's) Bivariate Gamma (pages
432--435)}\label{sec:cherrama}  Recall the Gamma distribution with
parameter $\theta>0$ is given by \be \twocase{p_Y(y) \ = \
}{\Gamma(\theta)^{-1} y^{\theta-1} e^{-y}}{if $y >
0$}{0}{otherwise.} \ee It has mean $\theta$, variance $\theta$,
and its moment generating function is \be M_Y(t) \ = \
(1-t)^{-\theta}, \ee which exists for all $t<1$. Let $Y_i$ for $i
\in \{0,1,2\}$ be independent Gamma distributed random variables
with parameters $\theta_i$, and for $i\in \{1,2\}$ set $X_i = Y_0
+ Y_i$. The density function of $(X_1,X_2)$ is \be
p_{X_1,X_2}(x_1,x_2) \ = \
\frac{e^{-(x_1+x_2)}}{\Gamma(\theta_0)\Gamma(\theta_1)\Gamma(\theta_2)}
\int_0^{\min(x_1,x_2)} y_0^{\theta_0-1}(x_1-y_0)^{\theta_1-1}
(x_2-y_0)^{\theta_2-1} e^{y_0} dy_0, \ee the bivariate gamma
density, equation (48.5). The correlation coefficient of $X_1$ and
$X_2$ is \be {\rm Corr}(X_1,X_2) \ = \
\frac{\theta_0}{\sqrt{(\theta_0+\theta_1)(\theta_0+\theta_2)}}.
\ee As $\theta_0 > 0$ (since $Y_0$ is Gamma distributed) the
correlation coefficient is positive, see (48.7). The moment
generating function is \be\label{eq:mgfbivar} M_{X_1,X_2}(t_1,t_2)
\ = \ (1-t_1-t_2)^{-\theta_0} (1-t_1)^{-\theta_1}
(1-t_2)^{-\theta_2} \ee and exists for all $(t_1,t_2)$ with
$t_1+t_2 < 1$ and $t_i < 1$, see (48.10).

\subsubsection{Multivariate Gamma
Distributions}\label{subsec:mvgdqazxsw}

We may generalize the arguments from Section \ref{sec:cherrama}
and consider the joint distribution of $X_p = \lambda_p(Y_0+Y_p)$
for $i \in \{1,\dots,P\}$ and $\lambda_p > 0$ with $Y_0, \dots,
Y_P$ independent Gamma distributed random variables with
parameters $\theta_0, \dots, \theta_P$. If $P=2$ Ghirtis has
called this the double-gamma distribution. For general $P$ it is
similar to Mathai and Moschopoulos' Multivariate Gamma
distribution (pages 465--470), and taking $\theta_p = 1$ we obtain
Freund's Multivariate Exponential distribution (pages 388--391).
The moment generating function is \bea\label{eq:GMMmultgammadistr}
M_{X_1,\dots,X_P}(t_1,\dots,t_P) & \ = \ & \E\left[e^{t_1 X_1 +
\cdots + t_P X_P}\right] \nonumber\\ & \ = \ & \E\left[ e^{\gl_1
t_1 (Y_0 + Y_1) + \cdots + \gl_P t_P (Y_0+Y_P)}\right] \nonumber\\
&=&\E\left[e^{(\gl_1 t_1 + \cdots + \gl_P t_P)Y_0}\right] \cdot
\E\left[e^{\gl_1 t_1 Y_1}\right] \cdots \E\left[e^{\gl_P t_P
Y_P}\right] \nonumber\\ &=& (1 - \gl_1 t_1 - \cdots - \gl_P
t_P)^{-\theta_0} (1 - \gl_1 t_1)^{-\theta_1} \cdots (1 - \gl_P
t_P)^{-\theta_P}, \eea which exists for all $(t_1,\dots,t_P)$ such
that $\gl_1 t_1 + \cdots + \gl_P t_P < 1$ and each $t_p <
\gl_p^{-1}$. For our applications such restrictions are harmless,
as in Theorem \ref{thm:MGFmultivariate} we evaluate the moment
generating function at $(-K_{i,1},\dots,-K_{i,P})$ and each
$K_{i,p} \ge 0$.

In fact, we may generalize even further.

\begin{lem}[(Generalized) Multivariate Gamma
Distribution]\label{lem:sjmillermultgamma} Let $Y_{0,1}, \dots,
Y_{0,M}, Y_1, \dots, Y_P$ be independent Gamma distributions with
parameters
$\theta_{0,1},\dots,\theta_{0,M},\theta_1,\dots,\theta_P$. For
$\lambda_{p,m}, \lambda_p \ge 0$ let \be X_p \ \ = \ \
\left(\gl_{p,1} Y_{0,1} + \cdots +\gl_{p,M} Y_{0,M}\right) \ + \
\gl_p Y_p, \ \ p \in \{1,\dots,P\}. \ee Then the moment generating
function is \bea\label{eq:GMMmultgammadistrG}
M_{X_1,\dots,X_P}(t_1,\dots,t_P) & \ = \ & \prod_{m=1}^M \left(1 -
\sum_{p=1}^P \gl_{p,m} t_p\right)^{-\theta_{0,m}}\ \cdot\
\prod_{p=1}^P (1 - \gl_p t_p)^{-\theta_p} \eea and exists for all
tuples $(t_1,\dots,t_P)$ where $\sum_{p=1}^P \gl_{p,m} t_p < 1$
for each $m$ and $t_p < \gl_p^{-1}$ for each $p$. For $r\neq s$
the covariances are  \bea {\rm Covar}(X_r,X_s) & \ = \
&\sum_{m=1}^M \gl_{r,m} \gl_{s,m} \theta_{0,m}, \eea and the
correlation coefficients are \be {\rm Corr}(X_r,X_s) \ = \ \frac{
\sum_{m=1}^M \gl_{r,m} \gl_{s,m} \theta_{0,m}}{\sqrt{\gl_{r,1}^2
\theta_{0,1} + \cdots + \gl_{r,M}^2 \theta_{0,M} + \gl_r^2 Y_r^2}
\ \sqrt{\gl_{s,1}^2 \theta_{0,1} + \cdots + \gl_{s,M}^2
\theta_{0,M} + \gl_s^2 Y_s^2}}.\ee
\end{lem}

\begin{proof} The moment generating function is \bea
M_{X_1,\dots,X_P}(t_1,\dots,t_P) & \ = \ & \E\left[e^{t_1 X_1 +
\cdots + t_P X_P}\right] \nonumber\\ & \ = \ & \E\left[ e^{
\sum_{p=1}^P (\gl_{p,1} Y_{0,1} + \cdots + \gl_{p,M} Y_{0,M} +
\gl_p Y_p) t_p  }\right] \nonumber\\
&=&\E\left[e^{ \sum_{p=1}^P \gl_{p,1} t_p Y_{0,1} }\right] \cdots
\E\left[e^{ \sum_{p=1}^P \gl_{p,M} t_p Y_{0,M}}\right] \cdot
\E\left[e^{\gl_1 t_1
Y_1}\right] \cdots \E\left[e^{\gl_P t_P Y_P}\right] \nonumber\\
&=& \prod_{m=1}^M \left(1 - \sum_{p=1}^P \gl_{p,m}
t_p\right)^{-\theta_{0,m}}\ \cdot\  \prod_{p=1}^P (1 - \gl_p
t_p)^{-\theta_p}, \eea which exists for tuples $(t_1,\dots,t_P)$
where $\sum_{p=1}^P \gl_{p,m} t_p < 1$ for each $m$ and $t_p <
\gl_p^{-1}$ for each $p$. For our applications such restrictions
are harmless, as in Theorem \ref{thm:MGFmultivariate} we evaluate
the moment generating function at $(-K_{i,1},\dots,-K_{i,P})$ and
each $K_{i,p} \ge 0$. The covariances and correlation coefficients
are easily determined in this case. As the $Y_{0,m}$ and $Y_p$ are
independent, for $r\neq s$ \bea {\rm Covar}(X_r,X_s) & \ = \ &
\E\left[(\gl_{r,1} Y_{0,1} + \cdots +\gl_{r,M} Y_{0,M} + \gl_r
Y_r)(\gl_{s,1} Y_{0,1} + \cdots +\gl_{s,M} Y_{0,M} + \gl_s
Y_s)\right] \nonumber\\ & = & \ - \ \E\left[\gl_{r,1} Y_{0,1} +
\cdots +\gl_{r,M} Y_{0,M} + \gl_r Y_r\right]\cdot \E\left[
\gl_{s,1}
Y_{0,1} + \cdots +\gl_{s,M} Y_{0,M} + \gl_s Y_s\right] \nonumber\\
& = & \E\left[(\gl_{r,1} Y_{0,1} + \cdots +\gl_{r,M} Y_{0,M})(
\gl_{s,1} Y_{0,1} + \cdots +\gl_{s,M} Y_{0,M})\right] \nonumber\\
& = & \ - \ \E\left[\gl_{r,1} Y_{0,1} + \cdots +\gl_{r,M} Y_{0,M}
\right]\cdot \E\left[ \gl_{s,1}
Y_{0,1} + \cdots +\gl_{s,M} Y_{0,M} \right] \nonumber\\
& = &  \sum_{u=1}^M \sum_{v=1}^M \gl_{r,u}\gl_{s,v}
\left(\E[Y_{0,u} Y_{0,v}] - \E[Y_{0,u}] \cdot \E[Y_{0,v}]\right)
\nonumber\\ &=&\sum_{m=1}^M \gl_{r,m} \gl_{s,m}{\rm Var}(Y_{0,m})
\nonumber\\ &=& \sum_{m=1}^M \gl_{r,m} \gl_{s,m} \theta_{0,m} \eea
and the correlation coefficient follows immediately. \end{proof}

Finally, to obtain an even more flexible distribution, we may
consider linear combinations of multivariate Gamma functions. The
methods of Section \ref{subsecgeneralize} are immediately
applicable and yield an extension of Theorem
\ref{thm:MGFmultivariate}.

\subsection{Computational Savings for Multivariate Distributions}

For a general multivariate distribution $G$ as in
\S\ref{subsec:generalG}, the efficiency of our expansion is
related to the rate of growth of the moments, which determines the
number of terms needed in the series expansions. However, if $G$
has a good closed-form expansion for its moment generating
function (as in \S\ref{subsec:goodmgfG}), then substantial
computational savings exist. We study the computational savings
for such $G$ below; we may take $G$ to be the bivariate gamma
distribution with MGF given by \eqref{eq:mgfbivar}, or the
multivariate generalizations of \eqref{eq:GMMmultgammadistr} or
\eqref{eq:GMMmultgammadistrG}.

Our assumptions on the moment generating function of $G$ imply the
conditions for Theorem \ref{thm:MGFmultivariate} are satisfied.
Thus we obtain a closed-form series expansion for $H_i$: \be H_i \
= \ \sum_{k_{i11}=0}^\infty \cdots \sum_{k_{iJN_i}=0}^\infty
(-1)^{\oa{k}\cdot\oa{1}} M_{\gbi}(-K_{i,1},\dots,-K_{i,P}), \ee
where as always \be  K_{i,p} \ = \ \sum_{j=1}^J \sum_{t=1}^{N_i}
(y_{ijt} + k_{ijt})x_{ijt,p}, \ \ \ \oa{k}\cdot\oa{1} \ = \
\sum_{j=1}^J \sum_{t=1}^{N_i} k_{ijt}. \ee Note again that $H_i$
depends weakly on $\oa{k}$; all that matters are the dot products
$\oa{k}\cdot \oa{x}_{i,p}$ and the parity of
$(-1)^{\oa{k}\cdot\oa{1}}$. We argue as in Theorem
\ref{thm:effunivargamma}. Let $r = (r_1,\dots,r_P)$ be a $P$-tuple
of non-negative integers. For each $\oa{k}$ we count the number of
solutions to the system of Diophantine equations
$\oa{k}\cdot\oa{x}_{i,p} = r_p$ ($p \in \{1,\dots,P\}$) while
recording the sign of $(-1)^{\sum_j\sum_t k_{ijt}} =
(-1)^{\oa{k}\cdot\oa{1}}$. Then we have

\begin{thm}\label{thm:effmultivargamma}
Set \begin{eqnarray}\label{eq:dpcountG} S(M) & \ = \ & \left\{v: v
=
(v_1,\dots,v_{M}), v_l \in \{0,1,2,\dots\} \right\} \nonumber\\
K_i(x,r,+) & \ = \ & \#\{k \in S(JN_i): \forall p \in
\{1,\dots,P\}, \oa{k}\cdot\oa{x}_{i,p} = r_p,
(-1)^{\oa{k}\cdot\oa{1}} = +1 \} \nonumber\\ K_i(x,r,-) & \ = \ &
\#\{k \in S(JN_i): \forall p \in \{1,\dots,P\},
\oa{k}\cdot\oa{x}_{i,p} = r_p, (-1)^{\oa{k}\cdot\oa{1}} = -1 \}.
\end{eqnarray}
Assume the $\beta_{i,p}$ are drawn from a one-sided multivariate
distribution with parameters $\oa{b}$ and moment generating
function $M_{\beta_i}(t_1,\dots,t_P)$ defined when each $t_p \le
0$. Then $P(Y|\Omega) = \prod_{i=1}^I \inti$ with
\begin{equation}\label{eqintigoodoneG}\inti \ = \ \sum_{r \in S(P)}
(K_i(x,r,+)-K_i(x,r,-)) \cdot M_{\gb_i}(-K_{i,1},\dots,-K_{i,P}),
\end{equation}
and the combinatorial and Diophantine estimates and bounds from
Appendix \ref{seccombdiophbounds} are still applicable, leading
again to enormous computational savings (after an initial one time
cost of determining the $K_i(x,r,\pm)$). \end{thm}

To gain additional savings in Theorem \ref{thm:effmultivargamma}
we may replace the multivariate distribution with a translated one
as in Section \ref{seccompissues}.

Further (at least if we use the multivariate distributions from
\S\ref{subsec:specialmultdistr}), $H_i$ is a sum of the moment
generating function, and the moment generating function is readily
differentiable in terms of its parameters. Thus we again obtain
closed-form expressions for the derivatives (see
\S\ref{seccompissues} and Appendix \ref{secnewtonprogram}), and
thus for certain data sets (where now the parameters may be
correlated) there is the possibility of using Newton's Method to
determine the optimal values.

Probably the most tractable and useful multivariate density will
be the multivariate gamma distribution from Lemma
\ref{lem:sjmillermultgamma}. While all covariances will be
non-negative, the moment generating function, covariances and
correlation coefficients are given by very simple formulas, and
are easily evaluated and easily differentiated. Moreover the
multivariate gamma distribution can take on a variety of shapes,
and as discussed in \S\ref{secgenfutureresearch} we may further
increase the admissible shapes by considering linear combinations
of multivariate gamma distributions.


\ \\
\begin{center}
\section{CONCLUSION}\label{conclusionsection}
\end{center}

In this research we obtain closed-form expansions for the
marginalization of the logit likelihood, allowing us to make
direct inferences about the population. In general these
expansions involve new special functions; however, in the case
where the distribution of heterogeneity follows a Gamma or
Multivariate Gamma distribution (or, in full generality, any
linear combination of multivariate distributions with a
closed-form moment generating function defined for all
non-positive inputs), by re-grouping the terms in the expansions
we obtain a rapidly converging series expansion of elementary
functions. We separate the calculations into two pieces. The first
piece is counting solutions to a system of Diophantine equations
(we are finding non-negative integer solutions $\oa{k}$ to
$\oa{k}\cdot\oa{x}_{i,p} = r_p$; these are linear equations with
integer coefficients); the second is evaluating certain
integrations, which depend only on $\Omega$ and the values of the
Diophantine sums.

The advantage of this approach is clear -- we need only do the
first calculations once. Thus, if we have $10^9$ or so operations
there, it is a one-time cost. When we need to evaluate the
functions at related points (say for the Newton's Method
maximization or at the grid points), we need only evaluate the
summations on $r = (r_1,\dots,r_P)$ in  \eqref{eqintigoodone},
\eqref{eqgoodYPexpansion} or \eqref{eqintigoodoneG}. This grouping
of terms is an enormous savings; we count the solutions to these
systems of equations once, and save the results as expansion
coefficients.

While this research has focused on one specific case, the logit
model, and two specific set of priors, the Gamma (if the response
coefficients are independent) and Multivariate Gamma (if there may
be correlations among the response coefficients) distributions,
our hope is that this research spurs others to consider deriving
closed-form solutions via expansions that can be made arbitrarily
close. In fact, closed-form expansions exist for any multivariate
distribution that has a closed-form moment generating function.
Thus our expansions can incorporate correlations among the
coefficients without sacrificing the computational gains.

As experience with pure simulation approaches shows, i.e. those
that are alternatives to that considered here, it is never a bad
thing to have an approach that can be used to explore the
parameter space (e.g. mode finding) in advance of running a
simulation routine. Whether it is to get good starting values, or
simply to understand the potentially multimodal nature of a
posterior surface, we hope that research such as this provides
value to researchers doing applied problems.

\appendix

\ \\
\begin{center}
\ \section{GAMMA FAMILY OF DISTRIBUTIONS}\label{ggamma}
\end{center}

Given the positivity restriction described in Section
\ref{geometric-series} for the $\beta_{i,p}$, we desired a family
of distributions defined on the positive real line that would be
extremely flexible, allowing for a variety of shapes of the
heterogeneity distribution; and, of course, be conjugate to the
geometric series expansion to the logit model.  The Generalized
Gamma family of distributions satisfies those requirements. As
this work concentrated on the Gamma distribution, we only describe
this case below.

\subsection{Plots}\label{subsec:plotsgammas}

We give a few plots of the  Gamma  distribution to illustrate the
richness of the family.
\begin{center}
\scalebox{.7}[.7]{\includegraphics{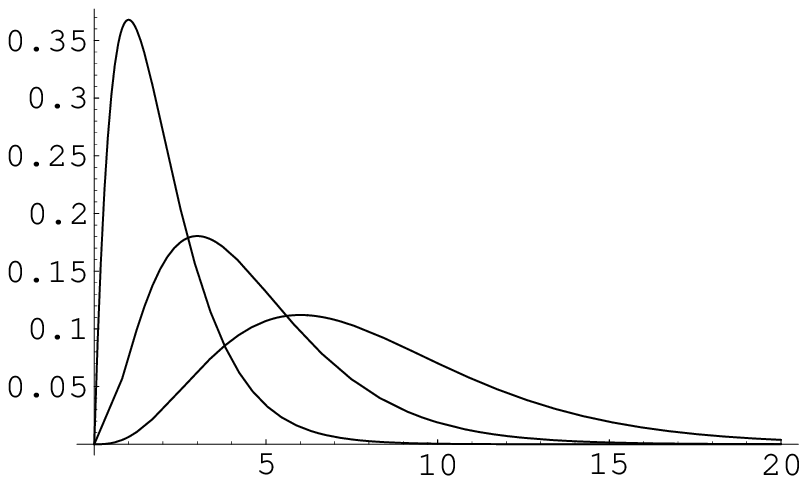}}\\
G(z;1,2), G(z;1.5,3), G(z;2,4).
\end{center}
While we develop the theory for $\beta_{i,p}$ drawn from a Gamma
distribution, we could use a weighted sum of Gamma distributions
as well.
\begin{center}
\scalebox{.75}[.75]{\includegraphics{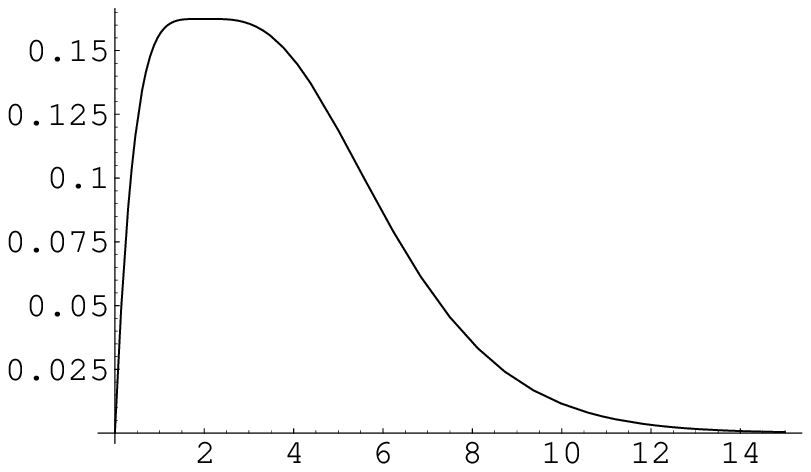}} \\ Sum of
Weighted Gamma distributions: $\frac{4}{10} \cdot G(z;1,2) +
\frac{6}{10} \cdot G(z;1,5)$:
\end{center}

\subsection{Integration Lemma}\label{subsectionggintproof}

We prove Lemma \ref{lemexpagainstggaisone}:

\begin{proof} We have
\begin{eqnarray} e^{-zd} G(z;b,n) & \ = \ & e^{-zd} \cdot \frac{1}{b
\Gamma(n)} \left(\frac{z}{b}\right)^{n-1} e^{-\frac{z}{b} }
\nonumber\\ & = & \frac{1}{b \Gamma(n)} \left(\frac{z}{b}
\right)^{n-1} e^{  - \frac{z}{ b / (1+bd) } } \nonumber\\ & = &
(1+bd)^{-n} \cdot \frac{1}{ \frac{b}{1+bd} \Gamma(n)} \cdot \left(
\frac{z}{ b / (1+bd) } \right)^{n-1} \cdot e^{ \frac{z}{b / (1+bd)
} } \nonumber\\ & = & (1+bd)^{-n} G\left(z, \frac{b}{1+bd},n
\right).  \end{eqnarray} As $b > 0$ and $d \ge 0$, $1 + bd > 0$,
and the above is well defined. Note $G(z; \frac{b}{1+bd},n)$ is
another Gamma distribution and therefore integrates to
$1$.\end{proof}


\ \\ \begin{center}
\section{MULTIVARIATE DENSITIES WITH CLOSED FORM MOMENT GENERATING FUNCTIONS}\label{sec:moremultdist}
\end{center}

In addition to the Multivariate Gamma distribution discussed in
detail in Section \ref{sec:covar}, we describe two additional
multivariate distributions that have closed-form expressions for
their moment generating functions. As such, these distributions
satisfy the conditions of Theorem \ref{thm:MGFmultivariate}, and
thus lead to closed-form series expansions. All page references
and equation numbers are from Kotz, Balakrishnan and Johnson 2000.
By no means is this list exhaustive, but rather representative of
those multivariate distributions which are  well suited to our
needs. Other distributions are Moran-Downton's Bivariate
Exponential (pages 371--377, especially (47.75) and (47.76)),
Freund's Multivariate Exponential (pages 388--391, especially
(47.85)), Kibble-Moran's Bivariate Gamma (pages 436--437),
Farlie-Gumble-Morgenstern Type Bivariate Gamma (pages 441--442,
especially (48.19) and (48.20)), and Mathai-Moschopoulos'
Multivariate Gamma (pages 465--470, especially (48.61) and
(48.62)). Other interesting distributions include truncated
multivariate normal distributions; however, as we require
one-sided distribution these are not as useful as those related to
the Gamma distributions.

\subsection{(Arnold and Strauss's) Bivariate Exponential (pages
370--371)}

Consider the joint probability density \be \twocase{
p_{X_1,X_2}(x_1,x_2) \ = \ }{A_{12} e^{-\gl_{12} x_1 x_2 - \gl_1
x_1 - \gl_2 x_2}}{if $x_1, x_2 > 0$}{0}{otherwise,} \ee where
$\gl_1, \gl_2, \gl_{12} > 0$ and $A_{12}$ is the normalization
constant. The moment generating function is given by \bea
M_{X_1,X_2}(t_1,t_2) & \ = \ & \E\left[e^{t_1X_1+t_2X_2}\right]
\nonumber\\ & = & A_{12}\int_0^\infty \int_0^\infty e^{-\gl_{12}
x_1 x_2 - \gl_1 x_1 - \gl_2 x_2 + t_1 x_1 + t_2 x_2} dx_1 dx_2
\nonumber\\ &=& A_{12}\int_0^\infty e^{-(\gl_2-t_2)x_2} \left[
\int_0^\infty e^{-(\gl_{12}x_2+\gl_1-t_1)x_1} dx_1 \right] dx_2
\nonumber\\ &=& A_{12} \int_0^\infty e^{-(\gl_2-t_2)x_2}
\frac{dx_2}{\gl_{12}x_2 + \gl_1 - t_1} \nonumber\\ &=&
\frac{A_{12}}{\gl_{12}} \int_0^\infty e^{-(\gl_2-t_2)x_2}
\frac{dx_2}{x_2 + (\gl_1-t_1)\gl_{12}^{-1}} \nonumber\\ &=&
\frac{A_{12}}{\gl_{12}} \int_0^\infty e^{-u}
\frac{du}{u+(\gl_1-t_1)(\gl_2-t_2)\gl_{12}^{-1}} \nonumber\\ &=&
\frac{A_{12}}{\gl_{12}}\ e^{-(\gl_1-t_1)(\gl_2-t_2)/\gl_{12}}\
{\rm Ei}\left(- \frac{(\gl_1-t_1)(\gl_2-t_2)}{\gl_{12}}\right),
\eea where \be {\rm Ei}(z) \ = \ - \int_{-z}^\infty e^{-t}
\frac{dt}{t} \ee is the exponential integral function (the
principal value is taken). The moment generating function exists
for $t_p < \gl_p$. For our applications such restrictions are
harmless, as in Theorem \ref{thm:MGFmultivariate} we evaluate the
moment generating function at $(-K_{i,1},\dots,-K_{i,P})$ and each
$K_{i,p} \ge 0$. The normalization constant can be determined by
setting $t_1 = t_2 = 0$: \be A_{12} \ = \ \gl_{12}
e^{\gl_1\gl_2/\gl_{12}}\ {\rm Ei}\left(-
\frac{\gl_1\gl_2}{\gl_{12}}\right)^{-1}. \ee

\subsection{(Freund's) Bivariate Exponential (pages 355--356)}
Freund considered the following situation: a two component
instrument has components with lifetimes having independent
density functions (when both are operating) of \be
\twocase{p_{X_p} \ = \ }{\alpha_p\ e^{-\ga_p x_p}}{if $x_p >
0$}{0}{otherwise,} \ee where $\ga_p > 0$; however, when one
component fails the parameter of the life distribution of the
other changes to $\alpha_k'$. Thus $X_1$ and $X_2$ are dependent
with joint density function \be \twocase{p_{X_1,X_2} \ = \
}{\ga_1\ga_2' e^{-\ga_2'x_2-\gamma_2 x_1}}{if $0 \le x_1 <
x_2$}{\ga_1'\ga_2 e^{-\ga_1'x_1-\gamma_1 x_2}}{if $0 \le x_2 <
x_1$,} \ee where $\gamma_p = \ga_1 + \ga_2 - \ga_p'$; see (47.25).
If $\gamma_p \neq 0$ then the marginal density of $X_p$ is \be
p_{X_p}(x_p) \ = \ \frac{(\ga_p-\ga_p')(\ga_1+\ga_2)}{\gamma_p}
e^{-(\ga_1+\ga_2)x_p} \ + \ \frac{\ga_p'\ga_{3-p}}{\gamma_p}
e^{-\ga_p' x_p}, \ \ \ x_p \ge 0. \ee As these are mixtures of
exponentials, this distribution is also called the bivariate
mixture exponential. The moment generating function is given by
\be M_{X_1,X_2}(t_1,t_2) \ = \ \frac{1}{\ga_1+\ga_2-t_1 -t_2}
\left(\frac{\ga_1'\ga_2}{\ga_1'-t_1} +
\frac{\ga_1\ga_2'}{\ga_2'-t_2}\right), \ee which converges for
$t_p < \ga_p'$ and $t_1+t_2 < \ga_1 + \ga_2$; see (47.28). For our
applications such restrictions are harmless, as in Theorem
\ref{thm:MGFmultivariate} we evaluate the moment generating
function at $(-K_{i,1},\dots,-K_{i,P})$ and each $K_{i,p} \ge 0$.
The correlation coefficient is given by \be {\rm corr}(X_1,X_2) \
= \ \frac{\ga_1'\ga_2' - \ga_1\ga_2}{\sqrt{(\ga_1'^2 + 2\ga_1\ga_2
+ \ga_2^2)(\ga_2'^2 + 2\ga_1\ga_2 + \ga_1^2)}} \ \in \
\left(-\frac13,\ 1\right), \ee see (47.31). Thus unlike many of
the other multivariate distributions, this model allows us to
study one-sided distributions with negative correlation.


\ \\ \begin{center}
\section{COMBINATORIAL AND DIOPHANTINE BOUNDS}\label{seccombdiophbounds}
\end{center}

We use the notation of Theorem \ref{thm:effunivargamma} and
Theorem \ref{thm:compefflincombuvg}:
\begin{eqnarray}\label{eq:dpcountapp} S(M) & \ = \ & \left\{v: v =
(v_1,\dots,v_{M}), v_l \in \{0, 1, 2, 3, \dots\} \right\}
\nonumber\\ K_i(x,r) & \ = \ & \#\{k \in S(JN_i): \forall p \in
\{1,\dots,P\}, \oa{k}\cdot\oa{x}_{i,p} = r_p\} \nonumber\\
K_i(x,r,+) & \ = \ & \#\{k \in S(JN_i): \forall p \in
\{1,\dots,P\}, \oa{k}\cdot\oa{x}_{i,p} = r_p,
(-1)^{\oa{k}\cdot\oa{1}} = +1 \} \nonumber\\ K_i(x,r,-) & \ = \ &
\#\{k \in S(JN_i): \forall p \in \{1,\dots,P\},
\oa{k}\cdot\oa{x}_{i,p} = r_p, (-1)^{\oa{k}\cdot\oa{1}} = -1 \},
\end{eqnarray} and let $K_i(r) = K_i(\oa{1},r)$.

For each $i$ we bound the number of solutions to $\oa{k}\cdot
\oa{x}_{i,p} = r_p$ for $p\in\{1,\dots,P\}$. Solutions to
Diophantine equations of this nature often crucially depend upon
the coefficients $x_{ijt,p}$. In expanding $P(Y|\beta)$ we can
trivially handle any terms with an $x_{ijt,p} = 0$. Thus, as we
assume $x_{ijt,p}$ is integral, in all arguments below we may
assume $x_{ijt,p} \ge 1$; if this assumption fails than trivial
book-keeping in our earlier expansions remove the sum over
$k_{ijt}$. The following result is immediate:

\begin{lem}\label{lem:simplercasexip1}
Let $x_{i,p} = (x_{i11,p},\dots,x_{iJN_i,p})$ be a $J\cdot N_i$
tuple of positive integers. Then $K_i(x,r,\pm) \le K_i(r)$.
\end{lem}

Thus by Lemma \ref{lem:simplercasexip1} instead of analyzing
$K_i(x,r,\pm 1)$ it suffices to bound the simpler $K_i(r)$.

\emph{For ease of exposition, we confine ourselves to the case
where the $\beta_{i,p}$ are drawn from a translated Gamma
distribution, $G(z-\epsilon;b_p,n_p)$, and we assume $x_{ijt,p}
\ge \delta$;} for example, we may take $\delta = 1$. Such bounds
do not exploit the cancellation in $K_i(x,r,+) - K_i(x,r,-)$
(though it is not unreasonable to expect square-root
cancellation). It is straightforward to generalize these arguments
to the Multivariate Gamma distribution (or linear combinations
thereof) from Lemma \ref{lem:sjmillermultgamma}.

Central in the arguments below are combinatorial results about
counting the number of representations of an integer as a sum of a
fixed number of integers. We briefly recall two useful results.

\begin{lem}\label{lemcookie} The number of ways to
write a non-negative integer $r$ as a sum of $P$ non-negative
integers is $\ncr{r+P-1}{P-1}$. \end{lem}

\begin{proof}[Sketch of the proof] Consider $r+P-1$ objects in a row.
Choosing $P-1$ objects partitions the remaining $r$ objects into
$P$ non-negative sets, and there are $\ncr{r+P-1}{P-1}$ ways to
choose $P-1$ objects from $r+P-1$ objects.
\end{proof}

\begin{lem}\label{lemcookiesum} The number of ways to
write a non-negative integer at most $R$ as a sum of $P$
non-negative integers is $\sum_{r=0}^R \ncr{r+P-1}{P-1} =
\ncr{R+P}{P}$.
\end{lem}

\begin{proof}[Sketch of the proof] Partition $R$ into $P+1$ sets
as in Lemma \ref{lemcookie}. As the last partition runs through
all numbers from $0$ to $R$ we get partitions of all numbers at
most $R$ into $P$ non-negative sets.
\end{proof}

To exploit the exponential decay in \eqref{eqgoodYPexpansion} from
the $\beta_{i,p}$ being drawn from translated Gamma distributions,
we must show that $K_i(r)$ does not grow too rapidly; we shall
show it grows at most polynomially in $r$. Note such arguments
ignore the decay of the $(1 + b_pY_{i,p} + b_pr_p)^{-n_p}$
factors. Assume we truncate our expansion by requiring $0 \le
k_{i11}+\cdots+k_{iJN_i} \le R$. As we assume that $x_{ijt,p} \ge
\delta$ and that we are using translated Gamma distributions, we
must bound
\begin{equation}\sum_{k_{i11}, \dots, k_{iJN_i}  \atop k_{i11} + \cdots +
k_{iJN_i} > R} \prod_{p=1}^P e^{- \gep \delta \sum_{j=1}^J
\sum_{t=1}^{N_i} k_{ijt} }. \end{equation} We use the notation of
Section \ref{seccompissues}. For any $r$, if each $\kijt \ge 0$,
then Lemmas \ref{lem:simplercasexip1} and \ref{lemcookie}
immediately yield

\begin{thm}\label{thm:polybound} We have
\begin{equation}K_i(r) \ = \ \#\{\oa{k}: k_{i11}+\cdots + k_{iJN_i} =
r\} \ = \ \ncr{r + JN_i - 1}{JN_i - 1}. \end{equation} Thus
$K_i(r) \le (r + JN_i-1)^{JN_i-1} / (JN_i-1)!$, which implies that
$K_i(r)$ grows at most polynomially. If $x_{ijt,p} \ge 1$ then
$K_i(x,r,\pm)$ grows at most polynomially.
\end{thm}

We conclude with some arguments and techniques that are specific
to having the exponential decay from the translated Gamma
distributions. These exploit improved bounds for summing $K_i(r)$
for $r$ in various ranges. We bound
\begin{equation}\label{eqsumedpr} \sum_{k_{i11}, \dots, k_{iJN_i}  \atop
k_{i11} + \cdots + k_{iJN_i} > R} \prod_{p=1}^P e^{- \gep \delta
\sum_{j=1}^J \sum_{t=1}^{N_i} k_{ijt} } \ = \ \sum_{r=R+1}^\infty
\ncr{r + JN_i - 1}{JN_i - 1} e^{-\gep \delta P r}. \end{equation}
By Lemma \ref{lemcookiesum} we have
\begin{eqnarray} \#\{\oa{k}: 0 \le k_{i11}+\cdots +
k_{iJN_i} \le  R\} & = &\ncr{R+JN_i}{JN_i}. \end{eqnarray}

\begin{rek} The number of $k$-tuples with $\sum_j \sum_t \kijt \le
R$ is $\ncr{R+JN_i}{JN_i}$. If we want the approximation from
looking at just these terms to be good, we need the sum in
\eqref{eqsumedpr} to be small. In this case, we initially need to
evaluate $\ncr{R+JN_i}{JN_i}$ terms, which leads to $R$ values to
store. In subsequent evaluations (note this encompasses not only
calculating $\inti$ but possibly also its partial derivatives
required for Newton's Method) we only have to read in $R$ values,
an enormous savings. The more varied the data $x_{ijt,p}$ is,
however, the more tuples of dot products must be stored. \end{rek}

To obtain a feel for these sizes, we tabulate the number of terms
arising from different values of $R, J$ and $N_i$. Note it is the
product $JN_i$ that matters, not the values of $J$ and $N_i$
separately.
\begin{center}\small
\begin{tabular}{c|c|c}
  $R$ & $JN_i$ & $\ncr{R+JN_i-1}{JN_i-1}$ \\
  \hline
        5  &     20  &    $10^{4.6}$ \\
        7  &     20  &    $10^{5.8}$ \\
        9  &     20  &    $10^{6.8}$ \\
        \hline
        5   &    30  &    $10^{5.4}$ \\
        7   &    30  &    $10^{6.9}$ \\
        9   &    30  &    $10^{8.2}$ \\
        \hline
        5   &    40  &    $10^{6.0}$ \\
        7   &    40  &    $10^{7.7}$ \\
        9   &    40  &    $10^{9.2}$ \\
\end{tabular}\normalsize
\end{center}
The largest term in the expansion of $\inti$ in Theorem
\ref{thm:univargamma} is when all $\kijt = 0$, giving $+1$. When
the $k$-sum is small (say of size $s$), we find terms of size
$e^{-\gep \delta P s}$. We have the following trivial estimate:
\begin{equation}\ncr{r + JN_i - 1}{JN_i - 1} \ \le \ (1 + 1)^{r+JN_i-1} \ = \
2^{JN_i-1} e^{r \log 2}. \end{equation} Assume $\gep \delta P >
\log 2$. Then the sum in  \eqref{eqsumedpr} is bounded by
\begin{eqnarray} \sum_{r=R+1}^\infty 2^{JN_i-1} e^{r \log 2} e^{-\gep \delta P
r} & \ \approx \ & 2^{JN_i-1} \int_{R}^\infty e^{-(\gep \delta P -
\log 2) r}dr \nonumber\\ & \approx & \frac{2^{JN_i-1} e^{-(\gep
\delta P - \log 2) \log R} }{\gep \delta P - \log 2}.
\end{eqnarray}
If $(\gep \delta P  - \log 2)\log R > JN_i \log 2$, the above is
small. Unfortunately, it might not be small compared to the
contributions from terms with a small $k$-sum (of size $s$); those
contribute on the order of $e^{-\gep \delta P s}$.

We perform a more delicate analysis by using dyadic decomposition,
breaking the sum over $r \ge R+1$ into blocks such as $2^m R \le r
\le 2^{m+1} R$, and using  Lemma \ref{lemcookiesum} in each block.
As the choice function $\ncr{r+JN_i-1}{JN_i-1}$ is monotonically
increasing in $r$, we find
\begin{eqnarray} \sum_{r=R+1}^\infty \ncr{r + JN_i - 1}{JN_i - 1} e^{-\gep
\delta P r} & \ < \ & \sum_{m=0}^\infty \sum_{r=2^m R}^{2^{m+1}R}
\ncr{r + JN_i - 1}{JN_i - 1} e^{-\gep \delta P r} \nonumber\\ & <
& \sum_{m=0}^\infty \ncr{2^{m+1}R + JN_i}{JN_i} e^{-\gep \delta P
2^m R} \nonumber\\ &  < &  \sum_{m=0}^\infty 2^{JN_i} e^{2^{m+1} R
\log 2} e^{-\gep \delta P 2^m R} \nonumber\\ & = & 2 e^{-((\gep
\delta P -2\log 2)R - JN_i\log 2)}. \end{eqnarray} This is small
if $(\gep \delta P - 2\log 2)R > JN_i \log 2$, allowing us to
replace the $\log R$ in $(\gep \delta P  - \log 2)\log R > JN_i
\log 2$ with $R$.

A slightly better savings is attainable by using instead
\begin{equation}\sum_{r=2^m R}^{2^{m+1}R} \ncr{r + JN_i - 1}{JN_i - 1} \ = \
\ncr{2^{m+1}R+JN_i}{JN_i} - \ncr{2^mR-1+JN_i}{JN_i}
\end{equation} and using polynomial (rather than exponential)
bounds. The main term is bounded by
\begin{equation}\twocase{\frac{(2^{m+1}R + JN_i)^{JN_i}}{(JN_i)!} \ < \
}{(2JN_i)^{JN_i}/(JN_i)!}{if $2^{m+1}R \le
JN_i$}{(2^{m+2}R)^{JN_i}/(JN_i)!}{if $2^{m+1}R > JN_i$}
\end{equation}

In order to deduce which of the many possible expansions is best,
and what size data sets are manageable, one needs to have explicit
values for $\gep, \delta$ and $P$; one can also try to exploit the
cancellation from the $(-1)^{\oa{k}\cdot\oa{1}}$ and the
denominator factors.

\ \\ \begin{center}
\section{APPLYING NEWTON'S METHOD TO THE
MARGINAL POSTERIOR}\label{secnewtonprogram}
\end{center}

Newton's Method yields a sequence of points $\vec{x}_k$ such that
$f(\vec{x}_k)$ converges to a local maximum of $f$. If $g_k$ and
$H_k$ are the gradient and Hessian of $f$ at $\oa{x}_k$, then
$\oa{x}_{k+1} = \oa{x}_k + \oa{p}_k$, where $\oa{p}_k$ satisfies
the \emph{linear} equation $H_k \oa{p}_k  = - \oa{g}_k$.

For our problem, $\Omega = (\oa{b},\oa{n})$. As the function we
want to maximize is a product of terms, we maximize $\log
f(\oa{b},\oa{n})$, as the logarithm converts the product in
\eqref{eqXXX} to a sum. To maximize
\begin{eqnarray} \log f(\oa{b},\oa{n}) & \ = \ & \log \prod_i
\inti(\oa{b},\oa{n}) \ = \ \sum_i \log \inti(\oa{b},\oa{n})
\end{eqnarray} we need the gradient and the Hessian as in standard applications
of Newton's method. The gradient is  \begin{eqnarray} \nabla \log
f(\oa{b},\oa{n}) \ = \ \frac{\nabla f(\oa{b},\oa{n})
}{f(\oa{b},\oa{n})} & = & \sum_i \frac{\nabla
\inti(\oa{b},\oa{n})}{\inti(\oa{b},\oa{n})},
\end{eqnarray} and the entries of the Hessian are
\be \pd{x} \grad \log f(\oa{b},\oa{n}) \ = \ \huge \sum_i
\normalsize \left[ \frac{ \pd{x} \grad
\inti(\oa{b},\oa{n})}{\inti(\oa{b},\oa{n})} - \frac{\grad \left[
\inti(\oa{b},\oa{n}) \right] \ \cdot \ \pd{x} \inti(\oa{b},\oa{n})
}{\inti^2(\oa{b},\oa{n})} \right], \ee where $\pd{x} = \pd{b_p}$
or $\pd{x} = \pd{n_p}$.

Straightforward differentiation gives the partial derivatives. The
advantage of using a Gamma distribution is the ease of
differentiating and evaluating these partials. We give exact,
infinite expansions; in practice, one truncates these expressions,
and the same Diophantine calculations and computational savings
for $\inti$ also hold for these derivatives.  Let
\begin{eqnarray}B(\oa{b},\oa{n},K(i)) & \ = \ & \prod_{p=1}^P
\left(1+b_pK_{i,p}\right)^{-n_p} \nonumber\\ \oa{k}\cdot \oa{1} &
= & k_{i11} + \cdots + k_{iJN_i}.
\end{eqnarray}

\begin{lem}[First Derivative Expansions]
\begin{eqnarray}\label{eqpartialsTinti} \frac{\partial \inti(\oa{b},\oa{n})}{\partial b_p} & \ = \
& - \sum_{k_{i11}=0}^\infty \cdots \sum_{k_{iJN_i}=0}^\infty
(-1)^{\oa{k}\cdot \oa{1}} \frac{K_{i,p} n_p}{1 + b_p K_{i,p}}
B(\oa{b},\oa{n},K(i)) \nonumber\\ \frac{\partial
\inti(\oa{b},\oa{n})}{\partial n_p} & \ = \ & -
\sum_{k_{i11}=0}^\infty \cdots \sum_{k_{iJN_i}=0}^\infty
(-1)^{\oa{k}\cdot \oa{1}} \log \left( 1 + b_p K_{i,p} \right)
\cdot
 B(\oa{b},\oa{n},K(i)). \end{eqnarray} \end{lem}

\noindent As $b_p, K_{i,p}$ are non-negative, the logarithms are
well defined above.

\begin{lem}[Second Derivative Expansions] In the expansions below,
$p \neq q$.
\begin{eqnarray} \frac{\partial^2 \inti(\oa{b},\oa{n})}{\partial
b_p^2} & = & \sum_{k_{i11}=0}^\infty \cdots
\sum_{k_{iJN_i}=0}^\infty (-1)^{\oa{k}\cdot \oa{1}}\frac{K_{i,p}^2
n_p(1 + n_p)}{(1 + b_p K_{i,p})^2} B(\oa{b},\oa{n},K(i)). \nonumber\\
\frac{\partial^2 \inti(\oa{b},\oa{n})}{\partial n_p^2} & = &
\sum_{k_{i11}=0}^\infty \cdots \sum_{k_{iJN_i}=0}^\infty
(-1)^{\oa{k}\cdot \oa{1}} \log^2 \left( 1 + b_p K_{i,p} \right)
\cdot  B(\oa{b},\oa{n},K(i)).\nonumber\\
 \frac{\partial^2 \inti}{\partial n_p \partial b_p} &=&
   \sum_{k_{i11}=0}^\infty \cdots
\sum_{k_{iJN_i}=0}^\infty (-1)^{\oa{k}\cdot \oa{1}} \Bigg[
\frac{K_{i,p}n_p}{1 + b_p K_{i,p}} \cdot \log (1 + b_p K_{i,p})
 \nonumber\\ & & \ \ \ \ \ \ \ \ \ \ \ \ \ \ \ \
 \ \ \ \ \ \ \ \ \ \ \ \ \ \ \ \ \ \ \ \ \ \ \ \ \ \ \  - \frac{K_{i,p}}{1 + b_p K_{i,p}}
\Bigg]\times B(\oa{b},\oa{n},K(i)) \nonumber\\ \frac{\partial^2
\inti(\oa{b},\oa{n})}{\partial b_p
\partial n_q} &  = & \sum_{k_{i11}=0}^\infty \cdots
\sum_{k_{iJN_i}=0}^\infty (-1)^{\oa{k}\cdot \oa{1}} \frac{K_{i,p}
n_p \log (1 + b_q K(i,q))}{1 + b_p K_{i,p}} B(\oa{b},\oa{n},K(i))
\nonumber\\ \frac{\partial^2 \inti(\oa{b},\oa{n})}{\partial b_p
\partial b_q} & \ = \ & \sum_{k_{i11}=0}^\infty \cdots
\sum_{k_{iJN_i}=0}^\infty (-1)^{\oa{k}\cdot \oa{1}} \frac{K_{i,p}
n_p}{1 + b_p K_{i,p}} \frac{K(i,q) n_q}{1 + b_q K(i,q)}
B(\oa{b},\oa{n},K(i)). \nonumber\\
 \frac{\partial^2 \inti}{\partial n_p \partial n_q} & = &
\sum_{k_{i11}=0}^\infty \cdots \sum_{k_{iJN_i}=0}^\infty
(-1)^{\oa{k}\cdot \oa{1}} \log \left( 1 + b_p K_{i,p} \right) \log
\left( 1 + b_q K(i,q ) \right)  B(\oa{b},\oa{n},K(i)).\ \ \ \ \ \
\end{eqnarray}
\end{lem}




\newpage

\begin{thebibliography}{PTW02} 



\bibitem{} J. Albert, Siddhartha Chib (1993), {\em Journal of the American
Statistical Association}, Vol. 88, No. 422., June, 669-679.


\bibitem{} R. Allen et al., Computational Science Education Project,
\emph{2.5.1 Newton Methods Overview},
\texttt{http://csep1.phy.ornl.gov/mo/node22.html}.

\bibitem{} G. Allenby and P.E. Rossi, (1993) \emph{A Bayesian Approach to Estimating Household Parameters},
Journal of Marketing Research, XXX, 171-182.

\bibitem{} G. M. Allenby, N. Aurora  and J.L. Ginter (1995), \emph{Incorporating
Prior Knowledge into the Analysis of Conjoint Studies}, Journal of
Marketing Research, 32, 152-162.

\bibitem{} E.T. Anderson and D.I. Simester (2001), \emph{Are Sale Signs Less Effective When More Products Have Them?}, Marketing Science, Vol. 20, No. 2, 121-142.

\bibitem{} B.L. Bayus (1992), \emph{Brand Loyalty and Marketing Strategy: An Application to Home Appliances}, Marketing Science, Vol. 11, No. 1, 21-38.

\bibitem{}M. Bierlaire, T. Lotan   and Ph. L. Toint (1997),  \emph{On the
overspecification of multinomial and nested logit models due to
alternative specific constants}, Transportation Science 31(4), 363--371.

\bibitem{} E.T. Bradlow, B.G.S. Hardie,  and P. Fader (2002),  \emph{Bayesian
Inference for the Negative-Binomial Distribution via Polynomial
Expansions}, Journal of Computational and Graphical Statistics,
Volume 11, Number 1, 189-201.

\bibitem{} P. E. Everson and E. T. Bradlow (2002), \emph{Bayesian Inference for
the Beta-Binomial Distribution via Polynomial Expansions}, Journal
of Computational and Graphical Statistics, Volume 11, Number 1, 202-207.

\bibitem{} A. E. Gelfand, S. Hills, A. Racine-Poon and A.F.M.
Smith (1990), \emph{Illustration of Bayesian Inference in Normal Data
Models Using Gibbs Sampling}, Journal Amer. Stat. Assoc., 85, 972-985.

\bibitem{} A. Gelman, J. Carlin, H.  Stern  and D. Rubin (1985), \emph{Bayesian Data Analysis},
CRC Press.

\bibitem{} A. Gelman and D. B. Rubin (1992), \emph{Inference from iterative
simulation using multiple sequences (with discussion)},
Statistical Science, 7, 457-511.


\bibitem{} J. Hausman  and D. McFadden (1984), {\em Specification Tests for the Multinomial Logit Model},
Econometrica, Vol. 52, No. 5., September, 1219-1240.

\bibitem{} W. Kamakura and G. Russell (1989), \emph{A Probabilistic
Choice Model for Market Segmentation and Elasticity Structuring,
Journal of Marketing Research}, November, 379-90.

\bibitem{} S. Kotz, N. Balakrishnan and N. L. Johnson,
\emph{Continuous Multivariate Distributions, Volume 1: Models and
Applications}, 2nd edition, John Wiley \& Sons, New York, 2000.

\bibitem{} R. MuCulloch and P. E. Rossi (1994), \emph{An Exact Likelihood Analysis of the Multinomial Probit
Model}, Journal of Econometrics, 64, 207-240.

\bibitem{} J. B. McDonald and R. J.  Butler (1990), \emph{Regression Models for
Positive Random Variables}, Journal of Econometrics, 41, 227-251.

\bibitem{} D. McFadden (1974), \emph{Conditional Logit Analysis of Qualitative Choice Behavior},
[PDF file, 3.2M] in P. Zarembka (ed.), FRONTIERS IN ECONOMETRICS, 105-142, Academic Press: New York.

\bibitem{} C. N. Morris (1983), \emph{Parametric Empirical Bayes Inference: Theory
and Applications}, Journal of the American Statistical
Association, Volume 78, Number 381, 47-65.

\bibitem{} D. G. Morrison an D. C. Schmittlein (1981),  {\em Predicting Future Random Events Based
on Past Performance}, Management Science, Vol. 27, No. 9., 1006-1023.

\bibitem{} D.  Revelt and K. Train (1998), \emph{Mixed Logit with Repeated Choices:
Households' Choices of Appliance Efficiency Level}, The Review of
Economics and Statistics, Vol. LXXX, No. 4, 647-657.

\bibitem{} G. J. Russell and R. Bolton (1988), \emph{Implications of Market Structure for Elasticity Structure}, Journal of Marketing Research, Vol. 25, No. 3, 229-241.

\bibitem{} D. C. Schmittlein, D. G. Morrison, and R. Colombo (1987),
{\em Predicting Future Random Events Based on Past Performance},
Management Science, Vol. 33, No. 1, 1-24.

\bibitem{} K. Train and D. McFadden (2000), {\em Mixed MNL Models for Discrete Response}, Journal of Applied
Econometrics, Vol. 15, No. 5, pp. 447-470.

\bibitem{} D. R. Wittink, M. J. Addona, W. J. Hawkes and J. C. Porter (1988),  \emph{SCAN*PRO(r): The Estimation,
Validation, and Use of Promotional Effects Based on Scanner Data}, February.

\end{thebibliography}
\end{document}